\documentstyle[12pt,ifthen]{article}
\newtheorem{theorem}{Theorem}[section]
\newtheorem{lemma}[theorem]{Lemma}
\newtheorem{proposition}[theorem]{Proposition}

\newtheorem{corollary}[theorem]{Corollary}
\newtheorem{remar}[theorem]{Remark}
\newenvironment{proof}{Proof:\ \ \ }{\QED}
\newenvironment{remark}{\begin{remar}\rm}{\end{remar}}

\newcommand{\QED}{{\unskip\nobreak\hfil\penalty50%
\hskip1em\hbox{}\nobreak\hfil $\Box$%
\parfillskip=0pt \finalhyphendemerits=0 \par\medskip\noindent}}
\newcommand{\bfind}[1]{\index{#1}{\bf #1}}

\newcommand{\sn}{\par\smallskip\noindent}

\newcommand{\bn}{\par\bigskip\noindent}
\newcommand{\pars}{\par\smallskip}

\newcommand{\sep}{^{\rm sep}}
\newcommand{\iT}{^{[i]}}

\newcommand{\chara}{\mbox{\rm char}\,}
\newcommand{\trdeg}{\mbox{\rm trdeg}\,}

\newcommand{\Spec}{{\mathrm{Spec}\,}}
\newcommand{\Frac}{{\mathrm{Frac}\,}}
\newcommand{\Gal}{{\mathrm{Gal}\,}}
\newcommand{\adresse}{\par\bigskip \small\rm
 Fraunhofer Institut Techno- und Wirtschaftsmathematik,\par
 Fraunhoferplatz 1, D--67663 Kaiserslautern,
 Germany\par
 email: knaf@itwm.fhg.de
  \par\bigskip
 Department of Mathematics and Statistics, 
 University of Saskatchewan, \par
 106 Wiggins Road, 
 Saskatoon, Saskatchewan, Canada S7N 5E6 \par
 email: fvk@math.usask.ca \ \ --- \ \ home page:
http://math.usask.ca/$\,\tilde{ }\,$fvk/index.html}
\font\tenlv=msbm10 scaled 1200
\font\sevenlv=msbm7 scaled 1200
\font\fivelv=msbm5 scaled 1200
\def\lv #1{{\mathchoice{{\hbox{\tenlv #1}}}{{\hbox{\tenlv #1}}}
{{\hbox{\sevenlv #1}}}{{\hbox{\fivelv #1}}}}}

\newcommand{\N}{\lv N}
\newcommand{\Q}{\lv Q}

\newcommand{\Z}{\lv Z}

\begin{document}
\title{Every place admits local uniformization in a finite
extension of the function field} 
\author{Hagen Knaf and Franz--Viktor Kuhlmann
\footnote{The second author thanks Peter Roquette, Mark
Spivakovsky, Bernard Teissier and Frans Oort for support and inspiring
conversations. Very special thanks to Hans Schoutens for many pleasant
and encouraging discussions.}}
\date{01.\ 11.\ 2006}
\maketitle
\begin{abstract}\noindent
We prove that every place $P$ of an algebraic function field $F|K$ of
arbitrary characteristic admits local uniformization in a finite
extension ${\cal F}$ of $F$. We show that ${\cal F}|F$ can be chosen
to be Galois, after a finite purely inseparable extension of the ground
field $K$. Instead of being Galois, the extension can also be chosen
such that the induced extension ${\cal F}P|FP$ of the residue fields is purely
inseparable and the value group of $F$ only gets divided by the
residue characteristic. If $F$ lies in the completion of an Abhyankar
place, then no extension of $F$ is needed. Our proofs are based solely
on valuation theoretical theorems, which are of particular importance in
positive characteristic. They are also applicable when working over a 
subring $R\subset K$ and yield similar results if $R$ is regular 
and of dimension smaller than $3$.
\end{abstract}
%
%
\section{Introduction and main results}
%
%
A place $P$ of an algebraic function field $F|K$ is said to admit local 
uniformization if there exists a $K$-variety $X$ having $F$ as its 
field of rational functions and such that the center $x\in X$ of $P$ 
on $X$ is a regular point. In [Z1], Zariski proved the Local Uniformization 
Theorem for places of algebraic function fields over base fields of 
characteristic 0. In [Z3], he uses this theorem to prove resolution of 
singularities for algebraic surfaces in characteristic 0, later on generalized 
to positive characteristic by Abhyankar [A1]. As the resolution 
of singularities for algebraic varieties of arbitrary dimension in positive 
characteristic is still an open problem, one is interested in generalizations 
of the Local Uniformization Theorem to positive characteristic. 
In this article we prove that every place of an algebraic function field of 
arbitrary characteristic admits local uniformization after a finite extension 
of the function field. This fact already follows from the results of de Jong 
[dJ] who proves resolution of singularities after a finite normal extension of 
the function field using results on moduli spaces of stable curves. However, 
we will give an entirely valuation theoretical proof which will provide 
important additional information about the finite extension used to achieve 
local uniformization. Our approach also applies to the case where the 
restriction of the place $P$ to $K$ is not the identity but is centered on a 
regular local ring $R\subset K$, $K=\Frac R$, of dimension 
$\dim R\leq 2$--thus including the arithmetic case of a discrete valuation 
ring $R$. In the latter case and for a function field $F|K$ of transcendence 
degree 1, Abhyankar [A2] has proved local uniformization (under some 
additional assumptions). If $R$ is a discrete valuation ring of a global 
field $K$ and for arbitrary transcendence degree of $F|K$, local uniformization 
after a finite extension of $F$ again follows from the results in [dJ].
\smallskip

Let $F|K$ be an algebraic function field equipped with a place $P$ whose 
restriction $P|_K$ to $K$ needs not be the identity. Local uniformization 
of $P$ is a statement about the valuation ring ${\cal O}_{P}$ associated with 
$P$. Accordingly throughout this article places $P$ and $P^\prime$ on the 
field $F$ inducing the same valuation ring are identified. By abuse of 
language the pair $(F|K,P)$ is called a \textbf{valued function field} keeping 
in mind the valuation $v$ of $F$ associated with $P$. The maximal ideal of 
the local ring ${\cal O}_{P}$ is denoted by ${\cal M}_{P}$ and the residue 
field of $P$ (or $v$) by $FP:={\cal O}_{P}/{\cal M}_{P}$.

Let $R\subseteq {\cal O}_{P}\cap K$ be a subring having field of fractions 
$\Frac R=K$. Given a 
separated, integral, finitely presented $R$-scheme $Y$ with $F=K(Y)$--an 
\textbf{$R$-model of $F|K$} for short--such that $P$ has center $y$ 
on $Y$ in the context of the resolution of singularities one searches for 
a birational morphism $X\rightarrow Y$ of $R$-models such that $P$ is 
centered in a regular point $x$ of $X$. Usually it is assumed that the schemes 
$X$ and $Y$ are noetherian, in the present article however we deal with the 
case of a non-noetherian valuation domain $R$ too. In that case one has to 
replace the requirement of being regular at the center $x\in X$ by 
smoothness of $X$ at the point $x$. For the 
valuation-theoretic approach presented in the sequel it is convenient to 
formulate the existence of the birational morphism $X\rightarrow Y$ in terms of 
the finite set of generators of the $R$-algebra ${\cal O}_Y(U)$ for a 
suitable open, affine neighborhood $U\subseteq Y$ of $y$. Doing so one 
arrives at the following notions: let $Z\subset {\cal O}_{P}$ be finite. 
The pair $(P,Z)$ is called \textbf{smoothly $R$-uniformizable} if there 
exists an $R$-model $X$ of $F|K$ such that $X\rightarrow\Spec R$ is smooth 
at the center $x\in X$ of $P$ on $X$ and $Z$ is contained in the local ring 
${\cal O}_{X,x}$ at $x$.
If $R$ is noetherian, the pair $(P,Z)$ is said to be 
\textbf{$R$-uniformizable} if $P$ is centered in a regular point $x\in X$ 
of an $R$-model $X$ of $F|K$ and  $Z\subset {\cal O}_{X,x}$ holds. The place 
$P$ is called (smoothly) $R$-uniformizable if the pair $(P,\emptyset )$ 
is (smoothly) $R$-uniformizable. The place $P$ is called \textbf{strongly 
(smoothly) $R$-uniformizable} if all pairs $(P,Z)$, $Z\subset{\cal O}_{P}$ 
finite, are (smoothly) $R$-uniformizable.
\smallskip

A natural approach to local uniformization is to consider stratifications 
of a valued function field $(F|K,P)$ essentially given through the choice 
of appropriate transcendence bases with respect to the place $P$: in general 
the inequality
\begin{equation}
\label{inequality}
\trdeg (FP|KP)+\dim (vF/vK\otimes_{\Z}\Q )\leq\trdeg (F|K)
\end{equation}
relates the transcendence degree $\trdeg (F|K)$ of $F|K$ with that of the 
residue field extension and with the rational rank of the abelian group 
$vF/vK$. The place $P$ is called an \textbf{Abhyankar place} if in 
(\ref{inequality}) equality holds. It is well-known that in every valued 
function field $(F|K,P)$ there exists an intermediate field $K\subset F_0
\subseteq F$ such that:
\begin{description}
\item[(S1)] the restriction $P|_{F_0}$ is an Abhyankar place of $F_0|K$ and 
$vF_0/vK$ is torsion-free,
\item[(S2)] the extension $FP|F_0P$ is algebraic and $vF/vF_0$ is a torsion 
group.
\end{description}
The field $F_0$ can be choosen to be a rational function field--see Theorem 
2.1 of [K--K] and Proposition \ref{tb} of the present article. Note also that if 
$P$ is not itself an Abhyankar place, then $F|F_0$ has positive transcendence 
degree.

In [K--K] valued function fields of the type appearing in (S1) are investigated: 
it is proved that an Abhyankar place $P_0$ of a function field 
$F_0|K$ is strongly $R$-uniformizable, where $R\subseteq K$ is a regular, 
local Nagata ring of Krull dimension $\dim R\leq 2$ dominated by ${\cal O}_{P_0}$, 
provided that the extension $FP_0|KP_0$ is separable and the valuation ring 
${\cal O}_{P_0}\cap K$ is defectless.

In the present article we study valued function fields $(E|K,P)$ as arising 
in (S2) with respect to smooth uniformizability of $P$ over the possibly 
non-noetherian valuation ring ${\cal O}_P\cap K$. We show that 
given a finite set $Z\subset {\cal O}_P$ there exists a finite extension 
${\cal E}|E$, a finite extension ${\cal K}|K$ within ${\cal E}$ and a 
prolongation ${\cal P}$ of $P$ to ${\cal E}$ such that the pair $({\cal P},Z)$ 
is smoothly $({\cal O}_{\cal P}\cap {\cal K})$-uniformizable.

The extension ${\cal E}|E$ can be choosen to be Galois. However for certain 
applications of local uniformization, e.g. to the model theory of fields in 
the spirit of [J--R] (cf.\ also [K5]), it is important to have a valuation 
theoretical control on the extension ${\cal E}|E$ and the residue field 
extension ${\cal E}{\cal P}|EP$ that we cannot obtain in the Galois case: we want 
to have ${\cal E}{\cal P}$ to be as close to $EP$ as possible, but in 
positive characteristic we may expect that we have to take a purely 
inseparable extension into the bargain. Therefore instead of choosing a 
suitable extension ${\cal E}|E$ within the separable closure $E\sep$ of $E$ 
we do the same within a separably tame hull of $E$: a valued field $(L,P)$ is 
called \bfind{separably tame} if it is henselian and its separable algebraic 
closure $L\sep$ equals the absolute ramification field of $(L,P)$. A 
\bfind{separably tame hull of the valued field $(E,P)$} is a field extension 
$E^{\rm st}|E$ equipped with an extension $P^{\rm st}$ of $P$ such that
$(E^{\rm st},P^{\rm st})$ is separably tame, 
$E^{\rm st}|E$ is separable-algebraic, 
$(v^{\rm st}E^{\rm st} /vE)$ is a $p$-group, 
and $E^{\rm st}P^{\rm st}|EP$ is a purely inseparable extension. Here $p$ 
denotes the characteristic of $EP$ respectively $p=1$ in the case of 
characteristic $0$. $v^{\rm st}$ is the valuation associated to the place 
$P^{\rm st}$. For basic properties of separably tame fields and the existence 
of separably tame hulls, see Subsection~\ref{subsectseptame}.
\begin{theorem}                         \label{main1}
Let $(E|K,P)$ be a separable, valued function field  such that $vE/vK$ is a
torsion group and $EP|KP$ is algebraic. Let $Z\subset {\cal O}_P$ be a finite
set. Let ${\cal P}$ be an extension of $P$ to the separable closure $E\sep$ 
of $E$. Then there exists a finite extension ${\cal E}|E$ within $E\sep$ and 
a finite extension ${\cal K}|K$ within ${\cal E}$ such that the function field 
${\cal E}|{\cal K}$ possesses an ${\cal O}_{\cal K}$-model $X$, 
${\cal O}_{\cal K}:={\cal O}_{\cal P}\cap {\cal K}$, with the properties:
\begin{itemize}
\item $X\rightarrow\Spec {\cal O}_{\cal K}$ is smooth at the center $x\in X$ 
of ${\cal P}$ on $X$,
\item every $z\in Z$ can be expressed as $z=uz^\prime$ with some $u\in
{\cal O}_{X,x}^{\times}$ and $z^\prime\in {\cal O}_{\cal K}$.
\end{itemize}
The extension ${\cal E}|E$ can be choosen to be either Galois or to be a 
subextension of a separably tame extension $E^{\rm st}|E$ within $E\sep$--for 
example a separably tame hull of $(E,P)$. If ${\cal E}|E$ is choosen to be 
Galois, then ${\cal K}|K$ can be choosen to be Galois too.

If $E_0|K$ is a subextension of $E|K$ such that $\trdeg E_0|K=\trdeg E|K -1$ 
and $E|E_0$ is separable, then ${\cal E}$ can be chosen to be a compositum 
$E.{\cal E}_0$, where ${\cal E}_0|E_0\,$ is a finite extension that is Galois 
respectively is contained in $E^{\rm st}$.
\end{theorem}
Let us return to a stratification $K\subset F_0\subset F$ satisfying the 
conditions (S1) and (S2) and assume in addition that $F|F_0$ is separable:
Theorem \ref{main1} yields finite extensions ${\cal F}|F$, ${\cal F}_0|F_0$ of a 
certain type such that $({\cal P}|_{\cal F},Z)$ is smoothly 
$({\cal O}_{\cal P}\cap {\cal F}_0)$-uniformizable for every finite set $Z\subset 
{\cal O}_P$. The place ${\cal P}|_{{\cal F}_0}$ is an Abhyankar place of the 
function field ${\cal F}_0|K$, thus the results on local uniformization 
of Abhyankar places obtained in [K--K] apply. Using a descend property of smooth 
algebras we can combine these facts to get general results about the 
uniformizability of pairs $(P,Z)$. Utilizing the statement in Theorem 
\ref{main1} concerning factorizations of the elements $z\in Z$ we can even extent 
these results to include monomialization of all $z\in Z$: let ${\cal O}$ be a 
commutative ring and $H\subseteq {\cal O}$. An element $a\in{\cal O}$ is called 
an \textbf{${\cal O}$-monomial in $H$} if
\[
a=u\prod\limits_{i=1}^d h_i^{\mu_i},\; u\in{\cal O}^\times ,\; h_i\in H,\;
\mu_i\in\N_0 ,\;  i=1,\ldots ,d,
\]
holds, where $\N_0:=\N\cup \{0\}$.

In the case of a valued function field $(F|K,P)$ with $P|_K={\rm id}_K$, in which 
case we also say that \textbf{$P$ is a place of $F|K$}, the combination of Theorem 
\ref{main1} with the results of [K--K] yields:
\begin{theorem}                                \label{main2}
Let $P$ be a place of the function field $F|K$ and 
let $Z\subset {\cal O}_P\,$ be a finite set. Let ${\cal P}$ be an extension of 
$P$ to the algebraic closure $\tilde{F}$ of $F$. Then there exist a finite 
purely inseparable extension ${\cal K}|K$ and a finite separable extension 
${\cal F}|F.{\cal K}$ such that the pair 
$({\cal P}|_{\cal F},Z)$ is smoothly ${\cal K}$-uniformizable.

More precisely let $F_0$ be an intermediate field of $F|K$ with the properties:
$P|_{F_0}$ is an Abhyankar place of $F_0|K$, $F|F_0$ is separabel, $vF/vF_0$ 
is a torsion group and $FP|F_0P$ is algebraic. Then there exist a finite 
purely inseparable extension ${\cal K}|K$, a finite separable extension 
${\cal F}|F.{\cal K}$, a finite extension ${\cal F}_0|F_0.{\cal K}$ within 
${\cal F}$ and a morphism $f:X\rightarrow X_0$ of ${\cal K}$-models of 
${\cal F}|{\cal K}$ and ${\cal F}_0|{\cal K}$ with the properties:
\begin{itemize}
\item $f$ is smooth at the center $x$ of ${\cal P}$ on $X$,
\item $X_0|{\cal K}$ is smooth at $f(x)$,
\item $\dim {\cal O}_{X,x}\geq\dim {\cal O}_{X_0,f(x)}=\dim (vF\otimes\Q )$,
\item all $z\in Z$ are ${\cal O}_{X,x}$-monomials in a regular parameter 
system of ${\cal O}_{X,x}$.
\end{itemize}
The extension ${\cal F}|F.{\cal K}$ can be choosen to be either Galois or to be 
a subextension of a given separably tame field $F^{\rm st}$ such that 
$F.{\cal K}\subseteq F^{\rm st}\subseteq (F.{\cal K})\sep$. In the first 
case the extension ${\cal F}_0|F_0.{\cal K}$ can be choosen to be Galois too.
\end{theorem}
The fact that in Theorem \ref{main2} one can choose $F^{\mathrm st}$ to be a 
separably tame hull of $F.{\cal K}$ has interesting consequences concerning 
the valuation theoretical control of ${\cal F}|F$ mentioned earlier:
\begin{corollary}                           \label{cor2}
If $\chara K=p>0$, then the valued extension $({\cal F}|F,{\cal P}|_{\cal F})$ 
can be choosen such that ${\cal F}{\cal P}|FP$ is a finite purely
inseparable extension and $v {\cal F}/vF$ is a finite $p$-group.
In particular one gets ${\cal F}{\cal P}=FP$ if $FP$ is perfect and 
$v{\cal F}=vF$ if $vF$ is $p$-divisible.

If $\chara K=0$, then one can take ${\cal F}$ to lie in the henselization
of $(F,P)$. In particular $v{\cal F}=vF$ and ${\cal F}{\cal P}=FP$ holds.
\end{corollary}
The assertion of the corollary in the case $p>0$ is a direct consequence 
of Theorem \ref{main2} and the definition of the separably tame hull. The 
assertion in the case $p=0$ can be considered as a weak version of Zariski's 
result on local uniformization [Z1]. It is a consequence of the fact that 
the henselization of $(F,P)$ is a separably tame hull--see Lemma \ref{tamerc0}.
\medskip

We turn to the case $P|_K\neq{\mathrm id}_K$, where we have to assume that 
the valued field $(K,P|_K)$ is defectless, that is that the fundamental equality 
of valuation theory holds in every finite extension $L|K$--see Section 
\ref{sectprel}. If the valuation associated to $P|_K$ is discrete, then 
defectlessness of $(K,P|_K)$ is equivalent to ${\cal O}_P\cap K$ being a Nagata 
ring. 
\begin{theorem}                                \label{main3}
Let $(F|K,P)$ be a valued function field such that $P|_K\neq {\mathrm id}_K$, 
$(K,P|_K)$ is defectless and $KP$ is perfect. Let ${\cal P}$ be an 
extension of $P$ to the separable closure $F\sep$ of $F$. 
Let $R\subseteq {\cal O}_{P}$ be a noetherian, regular local ring with maximal 
ideal $M={\cal M}_{P}\cap R$. Assume that $\Frac R=K$, $\dim R\leq 2$ and 
that $R$ is a Nagata ring if $\dim R=2$.

Then for every finite set $Z\subset {\cal O}_{P}$ there exists a finite 
separable extension ${\cal F}|F$ such that the pair $({\cal P}|_{\cal F},Z)$ is 
$R$-uniformizable.

More precisely let $F_0$ be an intermediate field of $F|K$ with the properties:
$P|_{F_0}$ is an Abhyankar place of $F_0|K$, $F|F_0$ is separabel, $vF_0 /vK$ 
is torsion-free,  $vF/vF_0$ is a torsion group and $FP|F_0P$ is algebraic. Then 
there exist a finite extension ${\cal F}|F$, a finite extension 
${\cal F}_0|F_0$ within ${\cal F}$ and a morphism $f:X\rightarrow X_0$ of 
$R$-models of ${\cal F}|K$ and ${\cal F}_0|K$ with the properties:
\begin{itemize}
\item $f$ is smooth at the center $x$ of ${\cal P}$ on $X$,
\item ${\cal O}_{X_0,f(x)}$ is regular,
\item $\dim {\cal O}_{X,x}\geq\dim {\cal O}_{X_0,f(x)}$, where
\[ 
\dim {\cal O}_{X_0,f(x)}=\left\{
\begin{array}{rl}
\dim (vF/vK\otimes\Q)+1 & \mbox{ if $\dim R=1$ or $\trdeg (KP|R/M)>0$}\\
\dim (vF/vK\otimes\Q)+2 & \mbox{ in the remaining cases}
\end{array}\right.
\]
\item all $z\in Z$ are ${\cal O}_{X,x}$-monomials in a regular parameter 
system of ${\cal O}_{X,x}$.
\end{itemize}
The extension ${\cal F}|F$ can be choosen to be either Galois or to be a 
subextension of a given separably tame field $F^{\rm st}$ such that 
$F\subseteq F^{\rm st}\subseteq F\sep$. 
In the first case the extension ${\cal F}_0|F_0$ can be choosen to be Galois too.
\end{theorem}
In [K--K] we have shown that Abhyankar places admit local uniformization 
without any extension of the function field. In [K5] a construction of places 
$P$ on a function field $F|K$ is given that yields non-Abhyankar places
which are still ``very close to'' Abhyankar places in the following sense: 
the valued field $(F,P)$ lies in the completion of a subfield $(F_0,P|_{F_0})$ 
such that $P|_{F_0}$ is an Abhyankar place. Therefore, it is important to know 
that also the latter places admit local uniformization without any extension 
of the function field. Here by ``completion'' we mean the completion with 
respect to the uniformity induced by the valuation: $(F,P)$ lies in 
the completion of $(F_0,P|_{F_0})$ if for every $a\in F$ and $\alpha\in vF$ 
there is some $b\in F_0$ such that $v(a-b)\geq\alpha$.
\begin{theorem}                                  \label{MTdens}
Let $(F|K,P)$ be a valued function field with the property that $(F,P)$ lies 
in the completion of a subfunction field $(F_0,P|_{F_0})$ such that $P|_{F_0}$ 
is an Abhyankar place of $F_0|K$, $vF_0 /vK$ is torsion-free and $F_0P|KP$ is 
separable.
\begin{enumerate}
\item If $P|_K={\mathrm id}_K$, then $P$ is strongly smoothly $K$-uniformizable 
and the conclusions of Theorem \ref{main2} concerning the existence and properties 
of the morphism $f: X\rightarrow X_0$ hold with ${\cal F}_0=F_0$ and ${\cal F}=F$.
\item Let $R\subset K$ be a subring of $K$ satisfying the requirements stated in 
Theorem \ref{main3}. If $P|_K\neq {\mathrm id}_K$ and $(K,P|_K)$ is defectless, 
then $P$ is strongly $R$-uniformizable and the conclusions of 
Theorem \ref{main3} concerning the existence and properties of the morphism 
$f: X\rightarrow X_0$ hold with ${\cal F}_0=F_0$ and ${\cal F}=F$.
\end{enumerate}
\end{theorem}
The results stated so far--in particular Theorem \ref{MTdens}--raise the question 
for necessary conditions for local uniformization without extending the function 
field. At least in the case of smooth uniformizability a condition in the same 
spirit as the major premise in Theorem \ref{MTdens} can be given: a valued 
function field $(F|K,P)$ is called 
\bfind{inertially generated} if it admits a transcendence basis $T$ such that 
$(F,P)$ lies in the absolute inertia field of $(K(T),P|_{K(T)})$. 
If it admits a transcendence basis $T$ such 
that $(F,P)$ lies in the henselization of $(K(T),P|_{K(T)})$, then we call it 
\bfind{henselian generated}.
\begin{theorem}                             \label{MT5}
Let $(F|K,P)$ be a valued function field such that $P$ is smoothly 
${\cal O}_K$-uniformizable. Then $(F|K,P)$ is inertially generated.
In particular $F|K$ and $FP|KP$ are separable. If in addition $FP=KP$, then 
$(F|K,P)$ is even henselian generated.
\end{theorem}
%
%
\section{Valuation theoretical preliminaries}         \label{sectprel}
In this section we review relevant facts from valuation theory in 
order to make the present article sufficiently self-contained. For basic 
facts from valuation theory we refer the reader to [EN], [R], [W] and [Z--S].
\subsection{Some fundamentals}
In the present article we formulate most of the results using the notion of a 
place of a field rather than that of a valuation to stress their geometric 
nature. It is well-known that the two notions essentially are synonymous to 
each other. Consequently by abuse of language we call a pair $(F,P)$ consisting 
of a field $F$ and a place $P$ of $F$ a \bfind{valued field}, keeping in mind 
the valuation associated to $P$, which we denote by $v\,$ or sometimes $v_P$ 
if explicit reference to the place $P$ is required. A \textbf{valued field 
extension} is a pair $(F|K,P)$, where $(F,P)$ is a valued field and $F$ is an 
extension field of $K$. The field $K$ is always understood to be equipped with 
the place $P|_{K}$, where we frequently suppress mentioning the restriction 
explicitely, that is we write $P$ instead of $P|_{K}$. 
If $F|K$ is finite respectively finitely generated, then we speak 
of a finite respectively finitely generated, valued field extension $(F|K,P)$.
The valuation ring of the valuation $v$ associated to $P$ is denoted by 
${\cal O}_P$ and its maximal ideal by ${\cal M}_P$. Additionally when 
considering intermediate fields $K\subseteq M\subseteq F$ of a valued field 
extension $(F|K,P)$ we use ${\cal O}_M :={\cal O}_P\cap M$ for the valuation 
ring of $v|_M$.

Throughout the article we identify places $P$ and $P^\prime$ of $F$ if they 
are inducing the same valuation ring of $F$. If that valuation ring is the 
field $F$ itself we call $P$ a \textbf{trivial place}. A trivial place is 
equivalent to the identity map of $F$. In particular if $(F|K,P)$ is a valued 
field extension such that $P|_K$ is an isomorphism of $K$, then we will assume 
that $P|_K={\mathrm id}_K$ and call $P$ a \textbf{place of $F|K$}.

Places operate on the right: the image of $f\in F$ under $P$ is denoted $fP$;
consequently $FP$ is the residue field ${\cal O}_{P}/{\cal M}_{P}$. The value 
group of the valuation $v$ associated to $P$ is denoted by $vF$ thus using the 
common convention $v(0)=\infty$.

For a valued extension $(L|K,P)$ the degree $f:=[LP:KP]$ is called 
\textbf{inertia degree} and $e:=(vL:vK)$ is the \textbf{ramification 
index}. If $L|K$ is finite, then $f$ and $e$ are finite too. More precisely 
if $P_1,\ldots ,P_g$ are the distinct extensions of $P|_K$ to $L$, then 
the {\bf fundamental inequality}
\begin{equation}                             \label{fundineq2}
[L:K]\>\geq\>\sum_{i=1}^g e_i f_i\; , 
\end{equation}
with $f_i=[LP_i:KP]$ and $e_i=(v_{P_i}L:vK)$, holds. 

A valued field $(K,P)$ is called 
\bfind{defectless} (or \bfind{stable}) if equality holds in (\ref{fundineq2}) 
for every finite extension $L|K$. As a consequence of the 
``Lemma of Ostrowski'' ([EN], [R]) a valued field with $\chara KP=0$ is 
defectless.

The effect of extending a place $P$ of a field $K$ to its separable 
closure $K\sep$ is described through the following fact: 
\begin{lemma}                              \label{KsacP}
Let $K$ be an arbitrary field and $P$ a non-trivial place on $K\sep$.
Then $v(K\sep )$ is the divisible hull $vK\otimes_{\Z}\Q$ of $vK$, and 
$K\sep P$ is the algebraic closure of $KP$.
\end{lemma}
For a proof see [K4], Lemma~2.16.

The valued extension $(L|K,P)$ is called \textbf{immediate} if $vL=vK$ and 
$LP=KP$.

A valued field $(K,P)$ is called \bfind{henselian} if it satisfies Hensel's 
Lemma; see [R] or [W]. The place $P$ then possesses a unique extension 
$P^\prime$ to every algebraic extension field $L$ of $K$ and $(L,P^\prime)$ 
is henselian too.

In general for every valued field $(K,P)$ there exists a henselian field 
$(K^h,P^h)$ and an embedding $i:K\rightarrow K^h$ such that $P=P^h\circ i$
with the following universal property: for every henselian extension 
$(L,P^\prime )$ of $(K,P)$ there exists a unique embedding 
$j:K^h\rightarrow L$ such that $P^h=P^\prime\circ j$. The valued field 
$(K^h,P^h)$ is uniquely determined up to a valuation-preserving 
$K$-isomorphism and is called the \textbf{henselization of $(K,P)$}. 
It can be contructed using ramification theory: define the 
\textbf{decomposition group} of an extension $P\sep$ of $P$ to $K\sep$ as 
$G_d:=\{\sigma\in\Gal (K\sep |K) : P\sep\circ \sigma =P\sep\}$. The fixed field 
of $G_d$ then is a henselization of $(K,P)$. The decomposition group contains 
the normal subgroup $G_i:=\{\sigma\in G_d : (\sigma (a)-a)P\sep=0\}$ called the 
\textbf{inertia group} of $P\sep$. The fixed field $K^i$ of $G_i$ equipped 
with the place $P^i:=P\sep |_{K^i}$ is henselian and is called the 
\textbf{absolute inertia field of $(K,P)$}; in the context of the present 
article the following property is relevant:
\begin{lemma}                               \label{ifcont}
Let $P$ be a place of $F|K$ and let $(F^i,P^i)$ denote the absolute inertia 
field of $(F,P)$. Then $K\sep\subset F^i$ holds. Further, if $FP|K$ is 
algebraic, then $(F.K\sep)P^i$ is the separable closure of $FP$.
\end{lemma}
\proof
By assumption $P|_K=\mathrm{id}_K$, hence $K\subseteq FP$. By general 
ramification theory we know  that $F^iP^i$ is separable-algebraically closed, 
thus $K\sep\subseteq F^iP^i$. Using 
Hensel's Lemma one can then construct a $K$-embedding $K\sep\hookrightarrow 
F^i$. Further $K\sep P^i\subseteq (F.K\sep)P^i$ and $K\sep\subseteq K\sep P^i$ 
by Lemma~\ref{KsacP}. As $F.K\sep|F$ is algebraic, so is $(F.K\sep)P^i|FP$. 
Therefore, if $FP|K$ is algebraic, then $(F.K\sep)P^i$ is algebraic over 
$K\sep$ and hence separable-algebraically closed. Since $(F.K\sep)P^i\subset
F^iP^i=(FP)\sep$, it follows that $(F.K\sep)P^i=(FP)\sep$.
\QED
%
%
\subsection{Transcendence bases of separable valued function fields}
The goal of the present section is to prove the existence of a transcendence 
basis of a valued function field $(F|K,P)$ that reflects basic properties 
of $P$ itself:
\begin{proposition}                               \label{tb}
Let $(F|K,P)$ be a valued function field and assume that $F|K$ is separable. 
Then there exists a separating transcendence basis of $F|K$ containing 
elements $x_1,\ldots,x_{\rho},y_1,\ldots,y_{\tau}$ such that:
\begin{itemize} 
\item The images of $vx_1,\ldots, vx_{\rho}$ under the natural map
$vF\rightarrow vF/vK\otimes\Q$ form a basis of the 
$\Q$-vector space on the right side,
\item $y_1P,\ldots,y_{\tau}P$ form a transcendence basis of $FP|KP$.
\end{itemize}
Here $v$ is the valuation of $F$ associated to $P$.
\end{proposition}
\begin{remark}
We do not know whether in addition to the assertion of the proposition, 
the $y_i$ can be chosen such that $y_1P,\ldots,y_{\tau}P$ form a separating
transcendence basis of $FP|KP$ if the latter extension is separable.
\end{remark}
To prove Proposition \ref{tb} we start with the case of a valued 
rational function field $(K(z)|K,P)$. The inequality (\ref{inequality}) then 
implies that the following three cases appear:
\begin{enumerate}
\item $K(z)P|KP$ is an algebraic extension and $vK(z)/vK$ is a torsion group,
\item $K(z)P|KP$ is a transcendental extension and $vK(z)/vK$ is a torsion group,
\item $K(z)P|KP$ is an algebraic extension and $vK(z)/vK$ is no torsion group.
\end{enumerate}
The first case can be characterized in terms of the behavior of the valued 
rational function field $(\tilde{K}(z)|\tilde{K},\tilde{P})$, where $\tilde{K}$ 
denotes the algebraic closure of $K$ and $\tilde{K}(z)$ is equipped with an 
arbitrary extension of the place $P$: using Lemma \ref{KsacP} we then see 
that the case 1 is equivalent to $(\tilde{K}(z)|\tilde{K},\tilde{P})$ being 
immediate. 
We use this fact in combination with the following easy to prove
\begin{lemma}                               \label{nomax}
The valued extension $(L|K,P)$ is immediate if and only if for every $z\in L$,
the set $\{v(z-a)\mid a\in K\}$ has no maximal element.
\end{lemma}
We then get--see [K4]: 
\begin{lemma}			\label{max}
For a valued rational function field $(K(z)|K,P)$ such that 
$(\tilde{K}(z)|\tilde{K},\tilde{P})$ is not immediate, there exists a monic 
irreducible polynomial $f\in K[X]$ that has a root in the set 
$\{a\in\tilde{K} : va=\max (v(z-b)\mid b\in \tilde{K})\}$. If $f$ has least 
degree among all such polynomials, then the following statements hold:
\begin{enumerate}
\item If $vK(z)/vK$ is no torsion group, then $vf(z)+vK$ is no torsion 
element.
\item If $K(z)P|KP$ is transcendental, then there is some $e\in\N$ and
some $d\in K$ such that $(df(z)^e)P$ is transcendental over $KP$.
\end{enumerate}
\end{lemma}
We deduce:
\begin{lemma}			\label{tb1}
In the situation of Lemma \ref{max} there exists a non-constant polynomial 
$h\in K[z]$ such that $K(z)|K(h)$ is separable and either $vh+vK$ is 
no torsion element in $vK(z)/vK$ or $hP$ is transcendental over $KP$.
\end{lemma}
\proof
If $vz$ is no torsion element of $vK(z)/vK$ or $zP$ is transcendental over 
$KP$, then $h:=z$ fullfills the requirements. Otherwise we treat the cases 
2 and 3 separately:

\textbf{case 3:} let $f\in K[X]$ be the irreducibel polynomial as defined in 
Lemma \ref{max}. We consider $h(X):=Xf(X)$: since 
$h^\prime (z)=zf^\prime (z)+f(z)\neq 0$ the element $z$ is a simple root of 
$h(X)-h(z)\in K(h(z))[X]$. Moreover $v(h(z))+vK=vz+v(f(z))+vK$; since  
$vz+vK$ is torsion by assumption, $h(z)$ fullfills all requirements by 
Lemma \ref{max}.

\textbf{case 2:} we consider $h(X):=Xdf(X)^e$, where 
$f\in K[X], d,e\in\N$ are choosen as in Lemma \ref{max}. Using the 
same argument as in the preceeding case we see that $z$ is a simple root 
of $h(X)-h(z)\in K(h(z))[X]$. Moreover $h(z)P=(zdf(z)^e)P$ is 
transcendental over $KP$ by Lemma \ref{max} and since $zP$ is assumed to 
be algebraic over $KP$.
\QED
\noindent
Proof of Proposition \ref{tb}: let $z_1,\ldots,z_n$ be a separating 
transcendence basis of $F|K$ and set $K_0:=K$, $K_i:=K(z_1,\ldots,z_i)$. 
Let $\tilde{P}$ be an extension of $P$ to the algebraic closure $\tilde{F}$ 
of $F$. For $i=1,...,n$ we consider the extension $(\tilde{K}_{i-1}(z_i)|
\tilde{K}_{i-1},\tilde{P})$, where $\tilde{K}_i$ denotes the algebraic 
closure of $K_i$ in $\tilde{F}$.

If $(\tilde{K}_{i-1}(z_i)|\tilde{K}_{i-1},P)$ is not immediate we choose 
$h_i\in K_{i-1}[z_i]$ as in Lemma \ref{tb1}. Otherwise, we set $h_i:=z_i$. 
The elements $h_1,\ldots,h_n$ then form a separating transcendence basis 
of $F|K$. The invariants $\rho:=\dim_{\Q}(vF/vK)\otimes\Q$ and 
$\tau:=\trdeg FP|KP$, where $v$ is the valuation associated to $P$, satisfy:
\[
\rho\,=\,\sum_{i=0}^{n-1}\dim_{\Q}(vK_{i+1}/vK_i\otimes\Q )
\mbox{ \ \ \ and \ \ \ }\tau\,=\, \sum_{i=0}^{n-1} \trdeg K_{i+1}P|K_iP\;.
\]
In view of the fact that
\[
\dim_{\Q}(vK_{i+1}/vK_i)\otimes\Q \,+\,\trdeg K_{i+1}P|K_iP\>\leq\> 
\trdeg K_{i+1}|K_i\>=\>1\;,
\]
we find that for precisely $\rho$ many values of $i$, $vh_i$ will be 
rationally independent modulo $v K_{i-1}\,$. Collecting all of these 
$h_i(z_i)$ and calling them $x_1, \ldots, x_{\rho}$ we thus obtain that 
$vx_1 ,\ldots, vx_{\rho}$ is a maximal set of elements in $vF$ rationally 
independent modulo $vK$.

Similarly, we find that for precisely $\tau$ many values of $\,i$,
the residues $h_iP$ will be transcendental over $K_{i-1}P$. Collecting all 
of these $h_i$ and calling them $y_1,\ldots,y_{\tau}$ we thus obtain that 
$y_1P,\ldots,y_{\tau}P$ form a transcendence basis of $FP|KP$.
\QED
%
%
\subsection{Separably tame fields}             \label{subsectseptame}
The \bfind{absolute ramification field} $K^r$ of a valued field $(K,P)$ is 
defined to be the fixed field of the group $G_r:=\{\sigma\in G_i : 
(\frac{\sigma (a)}{a}-1)P\sep=0\}$, where $P\sep$ is an extension of $P$ to 
the separable closure $K\sep$ of $K$. Let $p:=\max (1,\chara KP)$. By general 
ramification theory, $K\sep|K^r$ is a $p$-extension (cf.\ [EN]). Moreover 
$K^r$ contains the henselization $K^h$ of $(K,P)$ and is therefore henselian.
The valued field $(K,P)$ is called \bfind{separably tame} if it is henselian 
and satisfies $K\sep=K^r$--see also [K8].
\begin{lemma}                               \label{tamerc0}
Every henselian field of residue characteristic 0 is a separably tame
field.
\end{lemma}
Let $L|K$ be an algebraic extension. General ramification theory yields 
$L^r=L.K^r$ (cf.\ [EN]). Hence if $(K,P)$ is a separably tame field, then $L^r=L.K^r=L.K\sep=L\sep$. This proves:
\begin{lemma}				\label{algext}
Every algebraic extension of a separably tame field is separably tame.
\end{lemma}
A finite, separable extension $L$ of a separably tame field $(K,P)$ is a 
subextension of $K^r|K$ and thus, it satisfies the fundamental equality 
(cf.\ [EN]). This shows that every finite separable extension of a separably 
tame field is defectless. A valued field with this property is called a 
\bfind{separably defectless field}. So we note:
\begin{lemma}                               \label{tamehdp}
Every separably tame field is separably defectless.
\end{lemma}
The valued field $(K,P)$ is called \bfind{separable-algebraically maximal} 
if it admits no proper immediate separable-algebraic extension. Since the 
henselization is an immediate separable-algebraic extension (cf.\ [R]), 
we have:
\begin{lemma}                               \label{sam:hens}
Every separable-algebraically maximal field is henselian.
\end{lemma}
A finite, separable, immediate extension $L|K$ of a henselian, separably 
defectless field $(K,P)$ is trivial: $[L:K]= {\rm e}\cdot {\rm f}=1\cdot 1$.
Consequently:
\begin{lemma}                               \label{hdlam}
Every henselian, separably defectless field is separable-algebraically
maximal.
\end{lemma}
Combined with Lemma~\ref{tamehdp} this yields:
\begin{corollary}                           \label{corstam}
Every separably tame field is separable-algebraically maximal.
\end{corollary}
The subclass of separably tame fields within the class of 
separable-algebraically maximal fields can be characterized through 
conditions the value group and the residue field must satisfy--see 
also [K6]:
\begin{proposition}                    \label{tame}
Suppose that $P$ is a non-trivial place on $K$. Then $(K,P)$ is
separably tame if and only if it is separable-algebraically maximal,
$vK$ is $p$-divisible and $KP$ is perfect.
\end{proposition}
\proof
Assume first that $(K,P)$ is separably tame. By Corollary~\ref{corstam}, 
$(K,P)$ then is separable-algebraically maximal and by Lemma~\ref{KsacP}, 
$vK^r=vK\sep$ is divisible and $K^rP=K\sep P$ is algebraically closed, 
where $v$ denotes the valuation associated with the unique extension of 
$P$ to $K\sep$. This extension is denoted by $P$ again. General ramification 
theory tells us that every element of $vK^r/vK$ has order prime to $p$, and 
that $K^rP|KP$ is separable. Thus, $vK$ is $p$-divisible and $KP$ is perfect.

For the proof of the converse we start with the fact that for every 
henselian field $(K,P)$ there exists a subfield $W$ of $K\sep$ such that 
$W.K^r=K\sep$ and $W|K$ is linearly disjoint from $K^r|K$. This fact follows 
from Theorem~2.2 of [K--P--R] by Galois correspondence. 
Moreover Proposition~4.5 (ii) of [K--P--R] yields that $vW$ is the 
$p$-divisible hull of $vK$ and that $WP$ is the perfect hull of $KP$. 
In the present setting, as $vK$ is $p$-divisible and $KP$ is 
perfect, we conclude that $(W|K,P)$ is immediate. But as $(K,P)$ is
separable-algebraically maximal and $W|K$ is separable-algebraic, it follows 
that $W=K$. Consequently $K^r=W.K^r=K\sep$ showing that $(K,P)$ is separably 
tame.
\QED
\begin{corollary}                           \label{cortame}
If the valued field $(K,P)$ has $p$-divisible value group and perfect residue 
field, then every maximal immediate separable-algebraic extension of $(K,P)$ 
is a separably tame field. If $\chara KP= 0$ then already the henselization
$(K^h,P^h)$ is a separably tame field.
\end{corollary}
\proof
Let $(L,P^\prime )|(K,P)$ be a maximal immediate separable-algebraic 
extension. Then $(L,P^\prime )$ is separably-algebraically maximal, thus 
Proposition \ref{tame} yields the first assertion.

The henselization $(K^h,P^h )$ of $(K,P)$ is an immediate separable-algebraic 
extension. Lemmas~\ref{tamerc0} shows that $(K^h,P^h )$ is separably tame.
\QED
We next turn to the question under which conditions a subfield of a separably tame 
field inherits this property.
\begin{proposition}                                     \label{trac}
Let $(L,P)$ be a separably tame field. Assume that the subfield $K\subset L$ 
is separable-algebraically closed in $L$ and that $LP|KP$ 
is an algebraic extension, then $(K,P)$ is a separably tame field,  
the value group $vK$ is pure in $vL$ and $KP=LP$.
\end{proposition}
\proof
The field $K$ is separable-algebraically closed in the henselian field $L$, 
thus henselian too. Hensel's Lemma shows that $KP$ is 
separable-algebraically closed in $LP$. If $(K,P)$ is separably tame, then 
$KP$ is perfect by Lemma~\ref{tame}. Consequently we get $KP=LP$. In this 
situation Hensel's Lemma yields that $vL/vK$ is a $p$-group. On the other 
hand we know from Lemma~\ref{tame} that $vK$ is $p$-divisible. This shows
that $vK$ is pure in $vL$.

Altogether it remains to show that $(K,P)$ is separably tame. Considering Lemma 
\ref{tamerc0} from now on we can assume $p>0$. In order to prove that 
$K^r=K\sep$ holds, as in the proof of Proposition~\ref{tame} we choose a 
field $W\subseteq K\sep$ such that $W.K^r =K\sep$ and $W$, $K^r$ are linearly 
disjoint over $K$. We then have to show that $W=K$ holds.

Let $K'|K$ be a finite subextension of $W|K$. The degree of $K'|K$ is a power 
of $p$, since the Galois group $\Gal (K\sep |K^r )$ is known to be a $p$-group 
and $[K^\prime :K]=[K^\prime .K^r :K^r]$ by linear disjointness of $W$ and 
$K^r$ over $K$. 

The fields $K'$ and $L$ are linearly disjoint over $K$, since 
$K'|K$ is separable and $K$ is separable-algebraically closed in $L$.  Consequently 
$L':=L.K'$ satisfies $[L':L]=[K':K]$ and $K'$ is 
separable-algebraically closed in $L^\prime$. The extension $L'|L$ is 
separable and since $L$ is assumed to be separably tame we have $L\sep=L^r$.
We conclude $L'\subset L^r$. Since $L$ is henselian and the value group $v(L)$ is 
divisible by Propositon~\ref{tame}, general ramification theory yields that 
$L'P|LP$ is separable and that $[L^\prime :L]=[L^\prime P:LP]$ holds. 

Next utilizing Hensel's Lemma we see that $K'P$ is separable-algebraically 
closed in $L'P$. By hypothesis $LP|KP$ is an algebraic extension, therefore 
the same is true for $L'P|KP$ and thus for $L'P|K'P$. As a subextension of 
$WP|KP$ the extension $K'P|KP$ is purely inseparable. Now let $M|KP$ be a 
finite extension such that $M\subseteq LP$. Then
\[
[M.K'P:KP]_{\rm sep}^{}=[M.K'P:K'P]_{\rm sep}^{}[K'P:KP]_{\rm sep}^{}=1
\]
thus proving that $LP=L'P$. Consequently $L=L'$ and thus $K=K'$ as desired.
\QED
\textbf{Remark:} The preceeding proof is adapted from a proof that was 
given by F.~Pop for the case of tame fields.
\pars

A extension $(K^{\rm st}, P^{\rm st})$ of the valued field $(K,P)$ is called 
a \bfind{separably tame hull of} $(K,P)$ if it is a separably tame field with 
the following properties:
\begin{itemize}
\item $K^{\rm st}|K$ is separable-algebraic,
\item $v^{\rm st}K^{\rm st} /vK$ is a $p$-group,
\item $K^{\rm st}P^{\rm st}|KP$ is a purely inseparable extension.
\end{itemize}
These properties combined with Proposition~\ref{tame} imply that $v^{\rm st}
K^{\rm st}$ is the $p$-divisible hull of $vK$ and that $K^{\rm st}P^{\rm st}$ 
is the perfect hull of $KP$. 

A separably tame hull of a valued field $(K,P)$ always exists: in the case 
$p=1$ we can take the henselization $(K^h,P^h)$ of $(K,P)$. Otherwise let 
$W$ be an intermediate field of $K\sep|K^h$ such that $W$ and $K^r$ are 
linearly disjoint over $K^h$ and $W.K^r=K\sep$. Every maximal 
immediate separable-algebraic extension $K^{\rm st}$ of $W$ then is a 
separably tame hull of $(K,P)$ by Corollary \ref{cortame}. 
Unfortunately the separably tame hulls of $(K,P)$ 
are not unique up to valuation preserving isomorphism over $K$. However the 
failure of uniqueness does not matter for our use of separably tame hulls. 
%
%
\subsection{Kaplansky approximation}
For a polynomial $f\in K[z]$ in one variable over a field $K$ of arbitrary 
characteristic the $i$-th formal derivative $f\iT\in K[z]$ can be defined 
such that the following Taylor expansion holds (cf.\ [KA]):
\[
f(z)\>=\>f(a)+\sum_{i=1}^{\deg f} f\iT (a)(z-a)^i\;.
\]
Let $v$ be a valuation on $K(z)$. In this section we provide a result that 
allows to compute the value $vf(z)$ in terms of values derived from the 
Taylor polynomials after a suitable linear transformation of the variable $z$. 
Of course this is possible only if the valuation $v$ satisfies certain 
conditions; they were studied by Ostrowski and Kaplansky [KA]: 
let $(K(z)|K,P)$ be an immediate transcendental extension. 
The element $z$ induces the open sets $B(z,\alpha ):=\{a\in K :\; v(z-a)\geq
\alpha\}$, $\alpha\in vK$, in the uniform topological space $(K,v)$. 
Note that by the triangle inequality $B(z,\alpha )$ is a 
ball in $K$. These balls are interesting because of the particular behavior 
of maps $f:B(z,\alpha )\rightarrow K$ induced by polynomials $f\in K[z]$.
\begin{lemma}                               \label{close}
Let $(K(z)|K,P)$ be an immediate transcendental extension. Assume that
$(K,P)$ is a separable-algebraically maximal field or that $(K(z),P)$
lies in the completion of $(K,P)$. Then:
\begin{equation}	\label{trat}
\forall f\in K[z]\;\;\exists\alpha,\beta\in vK\;\;\forall a\in B(z,\beta ):\; 
vf(a)=\alpha.
\end{equation}
\end{lemma}
Kaplansky proved that if (\ref{trat}) does not hold, then there is a
proper immediate algebraic extension of $(K,P)$. If $(K(z),P)$ does not
lie in the completion of $(K,P)$, then using a variant of the Theorem on the 
Continuity of Roots this extension can be transformed into a proper immediate separable-algebraic extension (cf.\ [K6]). But such an extension cannot exist 
if we assume that $K$ be separable-algebraically maximal. 

If on the other hand $(K(z),P)$ lies in the completion of
$(K,P)$, then one can show that if $f$ does not satisfy (\ref{trat}),
then $vf(z)=\infty$. But this means that $f(z)=0$, contradicting
the assumption that $K(z)|K$ is transcendental.

\pars
We deduce the announced result about the computation of $vf(z)$:
\begin{lemma}                               \label{lvpol}
Let $(K(z)|K,P)$ be an immediate transcendental extension such that 
condition (\ref{trat}) holds. Then for every polynomial $f\in K[z]$
there exist $a,b\in K$ such that the values of the non-zero among the 
elements $f\iT (a)b^i$ are pairwise distinct and $v\tilde{z}=0$ for 
$\tilde{z}:= \frac{z-a}{b}$. In particular: 
\begin{equation}                            \label{vpol}
v f(z)\>=\>v(\sum_{i=0}^{\deg f} f\iT (a)(z-a)^i )\>=\>v
(\sum_{i=0}^{\deg f} f\iT (a)b^i \tilde{z}^i )\>=\>
\min_i\, (v (f\iT (a)b^i ))\;.
\end{equation}
If finitely many polynomials in $K[z]$ are given, then $a,b$ can be
chosen such that (\ref{vpol}) holds simultaneously for all of them.
\end{lemma}
\proof
Take finitely many polynomials $f_1,\ldots,f_n\in K[z]$. From
Lemma~\ref{close} we know that for the non-zero among the polynomials 
$f_j\iT$, $i,j\in\N$, there exist $\alpha_{ij},\beta \in vK$ such 
that: $\forall a\in B(z,\beta ):\; v f_j\iT (a)=\alpha_{ij}$ .
Since by Lemma~\ref{nomax} the set $\{v(z-a)\mid a\in K\}$ has no 
maximal element, we can choose $\beta\in vK$ so large that for 
$a\in B(z,\beta )$ and every fixed $j$, the values of the non-zero 
elements $f_j\iT (a)(z-a)^i$, $i\in\N$, are pairwise distinct. 

Having picked such
an element $a\in K$, we choose an element $b\in K$ such that $vb=
v(z-a)$. Then (\ref{vpol}) holds by the ultrametric triangle law.
\QED
%
%
\section{Smoothly uniformizable places}                \label{sectufe}
In this section we study valued function fields $(F|K,P)$ such that $(P,Z)$ 
is smoothly ${\cal O}_{K}$-uniformizable for some or all finite sets 
$Z\subset {\cal O}_{K}$, where ${\cal O}_{K}:={\cal O}_{P}\cap K$. We provide 
the basic properties of smooth uniformizability and prove a valuation-theoretic 
consequence: inertial generation of $F|K$--Theorem \ref{MT5}. Moreover we 
identify two classes of valued function fields whos members $(F|K,P)$ are 
strongly smoothly ${\cal O}_{K}$-uniformizable. One of these classes consists 
of the separable, immediate, valued function fields of transcendence degree 
one over a separably tame field $(K,P)$. The smooth uniformizability of the 
members of that class is a major building block of the proof of Theorem 
\ref{main1}.
\subsection{Basic properties}
Let $(F|K,P)$ be a valued function field. For the problem whether a pair 
$(P,Z)$ is smoothly $R$-uniformizable for some subring $R\subseteq 
{\cal O}_K$ it suffices to consider affine $R$-models $X=\Spec A$ of $F$, 
$A\subset F$ being a finitely presented $R$-algebra. Smoothness of $A$ at 
the center $q$ of $P$ on $X$ then means that there exists 
$f\in A\setminus q$ such that $A_f$ is $R$-flat and the rings 
$A_f\otimes_R \widetilde{k(p)}$, $\widetilde{k(p)}$ the algebraic closure of 
$k(p)=\Frac (R/p)$, are regular for all $p\in\Spec (R)$.
In the sequel we frequently need to contruct such an algebra $A$ within a 
given subring of $F$. The following structure theorem is 
particularly helpful in that respect--see [EGA IV], (17.11.4) for its proof:
\begin{theorem}                            \label{smoothness}
The $R$-algebra $A$ is smooth at $q\in\Spec A$ if and only if there exists
$u\in A\setminus q$ such that $A_u$ is an \'etale algebra over a polynomial
ring $R[x_1,\ldots ,x_d]$.
\end{theorem}
Recall that an $R$-algebra $A$ is called \textbf{standard-\'etale} if it admits 
a presentation of the form
\begin{equation}
\label{standard-etale}
0\rightarrow fR[X]_g\rightarrow R[X]_g\stackrel{\phi}{\rightarrow} A\rightarrow 0
\end{equation}
with $f,g\in R[X]$, $f$ monic and such that $\phi (f^\prime )\in A^\times$ for 
the derivative $f^\prime$ of $f$. Generalizing this definition we call an 
$R$-algebra $A$ \textbf{standard-smooth} if for some polynomial ring 
$T:=R[x_1,\ldots ,x_d]$ and some $h\in T$ the structure morphism $R\rightarrow A$ 
can be factored as
\[
R\rightarrow T_h\rightarrow A,
\] 
where $R\rightarrow T_h$ is the natural map and $T_h\rightarrow A$ is 
standard-\'etale. Consequently $A$ admits a presentation
\begin{equation}
\label{presentation}
0\rightarrow fT_h[X]_g\rightarrow T_h[X]_g\stackrel{\phi}{\rightarrow} 
A\rightarrow 0
\end{equation}
with $f,g\in T_h[X]$, $f$ monic, $\phi (f^\prime )\in A^\times$.
If $A$ is a domain, then the polynomial $f$ is prime in $T_h[X]_{g}$ but not 
necessarily in $T_h[X]$ itself. However if we assume $R$ to be normal, then 
$T_h$ is normal too, thus using Gau{\ss}' lemma we can choose $f$ to be prime 
in $T_h[X]$. 

Theorem \ref{smoothness} and the local structure theorem for \'etale algebras 
([Ray], Ch.~V, Thm.~1.) show that an $R$-algebra $A$ is smooth at 
$q\in\Spec A$ if and only if there exists some $u\in A\setminus q$ such that 
$A_u$ is standard-smooth.

Using standard-smooth algebras we prove that smoothness at a prime behaves 
well with respect to descent and ascent:
\begin{proposition}                         \label{descent 1}
Let $A|S$ be an extension of domains such that $S$ is normal and $A$ is a 
finitely presented $S$-algebra that is smooth at $q\in\Spec A$. Let 
$R\subseteq S$ be a subring of $S$ and let
$Z\subset A_q$ be a finite set. Then there exists a finitely generated ring 
extension $S_0 |R$ within $S$ with the property: for every normal
domain $S^\prime$ with $S_0\subseteq S^\prime\subseteq S$, there exists a 
finitely presented $S^\prime$-algebra $A^\prime\subseteq A_q$ that is smooth 
at $q^\prime :=qA_q\cap A^\prime$ and satisfies $Z\subset A^\prime_{q^\prime}$. 
Moreover for $F:=\Frac A$, $K:=\Frac S$ and $F^\prime :=\Frac A^\prime$ the 
relation $F=F^\prime .K$ holds. 
\end{proposition}
\proof
There exists $u\in A\setminus q$ such that $B:=A_u$ is a standard-smooth 
$S$-algebra. We choose a presentation of the form (\ref{presentation}) for 
$B|S$, where $T=S[x_1,\ldots x_d]$. Let $Z_1\subset S$ be the finite set 
of coefficients of $h\in T$ and of the coefficients $c\in T_h$ of $f$ and 
$g$. The condition $\phi (f^\prime )\in A^\times$ can be rewritten as
\begin{equation}
\label{unit-representation}
1=\phi (f^\prime\frac{t}{g^\ell}),\;\; t\in T_h[X], \ell\in\N ; 
\end{equation}
let $Z_2\subset S$ be the finite set of coefficients of the coefficients 
of $t$. Every $z\in Z$ can be expressed in the form
\begin{equation}
\label{z-representation}
z=\phi (\frac{p_z}{g^{k_z}})\phi (\frac{q_z}{g^{l_z}})^{-1},\;\; p_z,q_z\in T_h[X],\;k_z,l_z\in\N ,
\end{equation}
where $\phi (\frac{q_z}{g^{l_z}})\not\in qB$. 
Let $Z_3\subset S$ be the finite set of coefficients of the coefficients 
of the polynomials $\{p_z,q_z : z\in Z\}$.
Let $S_0:=R[Z_1\cup Z_2\cup Z_3]\subseteq S$ and consider a normal ring 
$S^\prime\subseteq S$ such that $S_0\subseteq S^\prime$. In the presentation 
(\ref{presentation}) choosen for $B|S$ we can then replace the ring 
$S$ by $S^\prime$ thus getting a presentation
\begin{equation}
0\rightarrow fT^\prime_h[X]_g\rightarrow T^\prime_h[X]_g\stackrel{\phi^\prime }
{\rightarrow} A^\prime\rightarrow 0
\end{equation}
with $T^\prime :=S^\prime [x_1,\ldots ,x_d]\subseteq T$. By construction 
$A^\prime$ is a standard-smooth $S^\prime$-algebra and the inclusion 
$T^\prime_h[X]_g\subseteq T_h[X]_g$ induces a homomorphism 
$A^\prime\rightarrow A$. We show that this map is injective: 
$T_h^\prime$ is normal since $S^\prime$ is so. An application of Gau{\ss}' 
lemma thus yields $fT_h[X]\cap T^\prime_h[X]  =fT^\prime_h[X]$ hence 
$fT_h[X]_g\cap T^\prime_h[X]_g  =fT^\prime_h[X]_g$.

By construction $Z\subset A^\prime_{q^\prime}$, $q^\prime :=A^\prime\cap A_q$, 
holds. Eventually with $K^\prime :=\Frac S^\prime$ we get 
$F^\prime .K=K^\prime (x_1,\ldots ,x_d,\phi^\prime (X)).K=
K(x_1,\ldots ,x_d,\phi (X))=F$ holds.
\QED
As for ascent we obtain:
\begin{proposition}                         \label{ascent 1}
Let $A|R$ be an extension of domains such that $R$ is normal and $A$ is a 
standard-smooth $R$-algebra. Let $S\supseteq R$ be a normal domain and 
assume that $F:=\Frac A$ and $L:=\Frac S$ are subfields of some field 
$\Omega$ such that $F$, $L$ are algebraically disjoint over $K:=\Frac R$. 
Then the compositum $A.S\subseteq F.L$ is a standard-smooth 
$S$-algebra. 
\end{proposition}
\proof
We choose a presentation of the form (\ref{presentation}) for $A|R$; 
it yields $A=T_h[x,g(x)^{-1}]$ with $T=R[x_1,\ldots ,x_d]$ a polynomial 
ring and $x=\phi (X)$, $g\in T_h[X]$. Consequently we 
get $A.S=T^\prime _h[x,g(x)^{-1}]$ with $T^\prime =S[x_1,\ldots ,x_d]$. 
The latter is a polynomial ring over $S$ since $F$ and $L$ are assumed to be 
algebraically disjoint over $K$. As mentioned earlier the normality of $T_h$ 
implies that the polynomial $f\in T_h[X]$ appearing in the presentation 
(\ref{presentation}) can be choosen to be the minimal polynomial of $x$ over 
$K(x_1,\ldots ,x_d)$. Let $f_1$ be the minimal 
polynomial of $x$ over $L(x_1,\ldots ,x_d)$. The normality of $T^\prime_h$ 
then yields $f_1\in T^\prime_h [X]$ and thus the exact sequence
\[
0\rightarrow f_1T^\prime_h[X]_g\rightarrow T^\prime_h[X]_g\rightarrow 
A.S\rightarrow 0.
\]
Moreover we have $f=f_1f_2$ for some $f_2\in T^\prime_h [X]$. Taking 
derivatives we obtain $f_1^\prime (x)f_2(x)=f^\prime (x)\in A^\times
\subseteq (A.S)^\times$, hence $f_1^\prime (x)\in (A.S)^\times$.
\QED
As an application of Proposition \ref{ascent 1} we can clarify some properties 
of smooth uniformizability over a valuation domain in a situation, where the 
constant field of the valued function field $(F|K,P)$ considered is extended:
\begin{proposition}                               \label{base change}
Let $(F|K,P)$ be a finitely generated, valued field extension. Let $L|K$ be a 
field extension and assume that $F$ and $L$ are subfields of some field 
$\Omega$ such that $F$ and $L$ are algebraically disjoint over $K$. 
Let ${\cal P}$ be an extension of $P$ to $F.L\subseteq\Omega$.
\begin{enumerate}
\item If $(P,Z)$ is smoothly ${\cal O}_K$-uniformizable, then $({\cal P},Z)$
is smoothly ${\cal O}_L$-uniformizable, where ${\cal O}_L:={\cal O}_{\cal P}
\cap L$.
\item Assume that $L|K$ is algebraic. If $({\cal P},Z)$ is smoothly 
${\cal O}_L$-uniformizable and $Z$ contains a set of generators of $F|K$, 
then there is a finitely generated subextension $M|K$ of $L|K$ such that 
$({\cal P}|_{F.M},Z)$ is smoothly ${\cal O}_M$-uniformizable. The field 
$M$ can be choosen to be algebraically closed in $F.M$.
\item Assume that $L|K$ is Galois. If $({\cal P},Z)$ is smoothly
${\cal O}_L$-uniformizable and $Z$ contains a set of generators of $F|K$, 
then there is a finite Galois subextension $N|K$ of $L|K$ containing $L\cap F$ 
such that $({\cal P}|_{F.N},Z)$ is smoothly 
${\cal O}_N$-uniformizable.
\end{enumerate}
\end{proposition}
\proof
\textbf{1.} There exists a standard-smooth ${\cal O}_K$-algebra $A\subset
{\cal O}_P$ such that $\Frac A=F$ and $Z\subset A_q$, $q:=A\cap {\cal M}_P$, 
hold. By Proposition \ref{ascent 1} the ${\cal O}_L$-algebra $B:=A.{\cal O}_L
\subset {\cal O}_{\cal P}$ is standard-smooth. It satisfies $\Frac B=F.L$.
Moreover for $q_B:=B\cap{\cal M}_{\cal P}$ the inclusion $Z\subset A_q
\subseteq B_{q_B}$ holds.

\textbf{2.} Let $A\subseteq {\cal O}_{\cal P}$ be a finitely presented 
${\cal O}_L$-algebra that is smooth at $q:={\cal M}_{\cal P}\cap A$ and 
satisfies $\Frac A=F.L$ and $Z\subset A_q$. Proposition \ref{descent 1} yields 
a finitely generated extension $S_0\subseteq {\cal O}_L$ of ${\cal O}_K$ such 
that for every valuation ring ${\cal O}_{M^\prime}\subseteq {\cal O}_L$ containing 
$S_0$ there exists a finitely presented ${\cal O}_{M^\prime}$-algebra $B\subseteq A_q$ 
that is smooth at $q_B:=qA_q\cap B$ and satisfies $Z\subset B_{q_B}$. We 
choose ${\cal O}_{M^\prime}$ such that $M^\prime=\Frac {\cal O}_{M^\prime}$ is a finitely generated 
extension of $K$. 
By assumption about $Z$ the field $E:=\Frac B$ contains $F$. Since $E|M^\prime$ is 
finitely generated, there exists a finitely generated extension $N|M^\prime$, 
$N\subseteq L$, such that $E\subseteq F.N=E.N$. By construction 
$({\cal P}|_E ,Z)$ is smoothly ${\cal O}_{M^\prime}$-uniformizable, hence applying (1) 
$({\cal P}|_{F.N},Z)$ is smoothly ${\cal O}_N$-uniformizable. For the  
algebraic closure $M$ of $N$ in $F.N$ we have $F.N=F.M$, thus 
we can apply (1) again to conclude that $({\cal P}|_{F.M},Z)$ is smoothly 
${\cal O}_{M}$-uniformizable.

\textbf{3.} Similarly to the first part of the proof of (2) we choose an 
${\cal O}_M$-algebra $B\subset A_q$ smooth at the center $q_B$ of ${\cal P}$ 
on $B$ and such that $Z\subset B_{q_B}$ holds. Since $(L\cap F)|K$ is finite 
we can assume that $M|K$ is a finite extension and $(L\cap F)\subseteq M$ holds. 
By assumption about $Z$ the 
field $E:=\Frac B$ contains $F$, thus the isomorphism of Galois groups
\[
\Gal (F.L|F)\rightarrow\Gal (L|L\cap F),\;\sigma\mapsto\sigma |_L
\]
yields $E=F.M^\prime$ for some finite extension $M^\prime|(L\cap F)$ such that 
$M^\prime\subseteq L$ and $M\subseteq M^\prime$. Let $N$ be the normal hull 
of $M^\prime |K$. Since $({\cal P}|_{E},Z)$ is 
smoothly ${\cal O}_{M}$-uniformizable by construction an application of (1) yields 
that $({\cal P}|_{E.N},Z)$ is smoothly ${\cal O}_{N}$-uniformizable. The equation 
$E.N=F.M^\prime .N=F.N$ concludes the proof.
\QED
Let $B$ be a smooth $R$-algebra and $C$ be a smooth $B$-algebra, then $C$ is 
a smooth $R$-algebra--[EGA IV],(17.3.3). Similarly if the $R$-algebra 
$B$ is regular at the prime $q_B$ and $q_C$ is a prime of the smooth 
$B$-algebra $C$ lying above $q_B$, then $C$ is regular at $q_C$--[EGA IV],(6.5.1). 
We next prove similar properties for (smooth) uniformizability.
\begin{proposition}                   \label{stacking}
Let $(F|L,P)$ be a finitely generated, valued field extension and assume 
that $(P,Z)$, $Z\subset{\cal O}_{P}$ finite, is smoothly 
${\cal O}_{L}$-uniformizable. Let $R$ be a subring of ${\cal O}_{L}$ and  
let $Z^\prime\subset {\cal O}_{L}$ be a finite set. Consider the following 
two cases:
\begin{description}
\item[case 1:] $P|_L$ is strongly smoothly $R$-uniformizable,
\item[case 2:] $R$ is noetherian and $P|_L$ is strongly $R$-uniformizable.
\end{description}
Then there exists a tower $R\subseteq B\subseteq C\subset {\cal O}_{P}$ of 
domains with fields of fractions $\Frac B=L$ and $\Frac C=F$ such that:
\begin{itemize}
\item the $R$-algebra $B$ is finitely presented in both cases, smooth in case 1,
and has the property that $B_{q_B}$ is regular in case 2,
\item the $B$-algebra $C$ is finitely presented  in both cases, smooth in case 1, 
and has the property that $C_{q_B}$ is a smooth $B_{q_B}$-algebra in case 2,
\item in both cases $Z^\prime\subset B_{q_B}$ and $Z\subset C_{q_C}$ hold.
\end{itemize}
Consequently the pair $(P,Z)$ is smoothly 
$R$-uniformizable in case 1 and $R$-uniformizable in case 2.
\end{proposition}
\proof
Take a standard-smooth ${\cal O}_{L}$-algebra $A\subset{\cal O}_{P}$ 
with the properties $\Frac A=F$ and $Z\subset A_{q_A}$, where 
$q_A:={\cal M}_{P}\cap A$. After choosing a presentation of $A|{\cal O}_{L}$ 
of the form (\ref{presentation}), we can define the finite sets 
$Z_1,Z_2,Z_3\subset {\cal O}_{L}$ as in the proof of Proposition 
\ref{descent 1}. By assumption there then exists a finitely presented 
$R$-algebra $B\subseteq{\cal O}_{L}$ with the properties $\Frac B=L$ and 
$Z^\prime\cup Z_1\cup Z_2\cup Z_3\subset B_{q_B}$, where $q_B:={\cal M}_{P}
\cap B$. Furthermore $B|R$ is smooth at $q_B$ in case 1 and $B_{q_B}$ is 
regular in case 2. 

In case 1 by passing from $B$ to a suitable localization $B_u$, $u\not\in q_B$, 
we can assume that $Z^\prime\cup Z_1\cup Z_2\cup Z_3\subset B$ and that $B$ is 
a smooth $R$-algebra. In particular $B$ is normal, hence as in the proof of 
Proposition \ref{descent 1} we can construct a standard-smooth $B$-algebra 
$C\subseteq A$ such that $\Frac C=\Frac A$ and $Z\subset C_{q_C}$ hold.
Since $C$ then is a smooth $R$-algebra the proposition is proved in case 1.

In case 2 since $B_{q_B}$ is normal as in the proof of 
Proposition \ref{descent 1} we can construct a standard-smooth $B_{q_B}$-algebra 
$C^\prime\subseteq A$ such that $\Frac C^\prime=\Frac A$ and $Z\subset 
C^\prime_{q_{C^\prime}}$ hold, where $q_{C^\prime}:={\cal M}_P\cap C^\prime$.
For $C^\prime =B_{q_B}[x_1,\ldots ,x_r]$ we set $C :=B[x_1,\ldots ,x_r]$; then
$C_{q_B}=C^\prime$ and $C_{q_C}=
C^\prime_{q_{C^\prime}}$. The smoothness of $C^\prime |B_{q_B}$ and 
the regularity of $B_{q_B}$ thus imply the regularity of $C_{q_C}$ and the 
proof is complete in the case 2.
\QED%
\begin{corollary}                               \label{trans}
Let $(F|K,P)$ be a finitely generated, valued field extension and $L$ an 
intermediate field of $F|K$. If $P|_L$ is strongly smoothly 
${\cal O}_{K}$-uniformizable and $P$ is strongly smoothly 
${\cal O}_{L}$-uniformizable, then $P$ is strongly smoothly 
${\cal O}_{K}$-uniformizable.
\end{corollary}
\subsection{Inertially generated function fields}
Applying the results of the previous subsection we are now able to provide 
the proof of Theorem \ref{MT5}.
\begin{lemma}                              \label{etale+smoothly uniform}
Let $(F|L,P)$ be a finite valued field extension and let $R\subseteq 
{\cal O}_{L}$ be a subring of $L$ with $\Frac R=L=\Frac {\cal O}_L$.
Then the following statements hold:
\begin{enumerate}
\item If $P$ is smoothly $R$-uniformizable, then ${\cal O}_{P}|
{\cal O}_{L}$ is local-\'etale.
\item $P$ is strongly smoothly ${\cal O}_{L}$-uniformizable if and only if 
${\cal O}_{P}|{\cal O}_{L}$ is local-\'etale.
\item ${\cal O}_{P}|{\cal O}_{L}$ is local-\'etale if and only if $(F,P)$ lies in 
the absolute inertia field of $(L,P)$.
\end{enumerate}
\end{lemma}
\proof
\textbf{1.:} Let $A\subset {\cal O}_{P}$ be a finitely presented $R$-algebra 
that is smooth at $q:={\cal M}_{P}\cap A$ and satisfies $\Frac A=F$. Since 
$F|L$ is algebraic Theorem \ref{smoothness} shows that we can assume $A|R$ to 
be standard-\'etale. An application of Proposition \ref{ascent 1} yields that 
the ${\cal O}_{L}$-algebra $B:=A.{\cal O}_{L}\subseteq {\cal O}_{P}$ is 
standard-\'etale too. Hence it suffices to prove $B_{q_B}={\cal O}_{P}$ for 
$q_B:={\cal M}_{P}\cap B$. Indeed as an \'etale extension of the normal domain 
${\cal O}_{L}$ the domain $B$ is normal too [Ray], Ch.~VII, Prop.~2. Hence 
$B_{q_B}$ is a normal local extension of ${\cal O}_{L}$ in the finite extension 
$F|L$ and thus a valuation domain contained in ${\cal O}_{P}$. Since the 
valuation rings of $F$ that are local extensions of ${\cal O}_{L}$ are 
pairwise incomparable with respect to inclusion we get $B_{q_B}={\cal O}_{P}$ 
as desired.

\textbf{2.:} The remaining implication $\Leftarrow$ is obvious.

\textbf{3.:}  See [Ray], Ch.~X., Thm.~1.
\QED%
\textbf{Proof of Theorem \ref{MT5}}: let $X=\Spec A$ be an affine 
${\cal O}_{K}$-model of the valued function field $(F|K,P)$, that is smooth 
at the center $q:={\cal M}_{P}\cap A$ of $P$ on $X$. By Theorem 
\ref{smoothness} we can assume that $A$ is an \'etale extension of a 
polynomial ring ${\cal O}_{K}[x_1,\ldots ,x_d]\subseteq A$. In particular 
the set $T:=\{x_1,\ldots ,x_d\}$ forms a transcendence basis of $F|K$ and 
$(P,\emptyset )$ is smoothly ${\cal O}_{K}[x_1,\ldots ,x_d]$-uniformizable. 
An application of Lemma \ref{etale+smoothly uniform}, (1) yields that 
${\cal O}_{P}|{\cal O}_{K(T)}$ is local-\'etale. Thus $F$ lies in the absolute 
inertia field of $(K(T),P)$ by (3) of the same lemma. 
Finally assume that $FP=KP$ holds. Then $F$ is an extension of $K(T)$
within the inertia field of $(K(T),P)$ such that $FP=K(T)P$. Thus $F$ 
must lie in the henselization of $(K(T),P)$.\QED
\subsection{Immediate function fields of transcendence degree one}
It is tempting to try to prove the reversed implication in Theorem \ref{MT5}. 
This however amounts to proving the ${\cal O}_K$-uniformizability of all 
valued, rational function fields $(K(T),P)$. In the sequel we present a case, 
where rational function fields of one variable are strongly smoothly 
${\cal O}_K$-uniformizable and draw some conclusions utilizing a 
structure theorem for immediate function fields over a separably tame field:
\begin{theorem}                \label{stt3}
Let $(F|K,P)$ be an immediate, valued function field of trans\-cendence 
degree 1 and assume that $(K,P)$ is separably tame. If $F|K$ is separable, 
then there exists $x\in F$ such that $(F,P)$ lies in the henselization 
$(K(x)^h,P^h)$, that is $(F|K,P)$ is henselian generated.
\end{theorem}
For the case $\chara KP=0$ the assertion is a direct consequence of the fact 
that every such field is defectless--in fact every $x\in F\setminus K$ will 
then do the job. In contrast to this, the case $\chara KP\neq 0$ requires a 
much deeper structure theory of immediate algebraic extensions of henselian 
fields, in order to find suitable elements $x$. For the proof of the theorem 
in this case see [K8].

\begin{lemma}                               \label{ist}
Let $(K(x)|K,P)$ be an immediate, transcendental extension possessing the 
property (\ref{trat}) stated in Lemma \ref{close}, then $P$ is strongly 
smoothly ${\cal O}_{K}$-uniformizable.
\end{lemma}
\proof
Let $z_1,\ldots,z_m\in {\cal O}_{P}\,$ and write $z_j=f_j(x)/g_j(x)$ with 
polynomials $f_j,g_j\in K[x]$. We apply Lemma~\ref{lvpol} to these finitely 
many polynomials and choose $\tilde{x}=\frac{x-a}{b}$, $a,b\in K$, according 
to this lemma. 
Then by (\ref{vpol}), for every $j$ we can find $i_j,k_j$ such that 
\[
vf_j(x)= v f_j^{[i_j]} (a)\,b^{i_j} =\min_i v f_j\iT (a)\,b^i
\mbox{ and }
vg_j(x)= vg_j^{[k_j]} (a)\,b^{k_j}=\min_i vg_j\iT (a)\,b^i.
\]
Thus we can write
\begin{equation}                   \label{represent}
z_j\>=\>\frac{f_j^{[i_j]}(a)\,b^{i_j}} {g_j^{[k_j]}(a)\,b^{k_j}}
\cdot \frac{\tilde{f}_j(\tilde{x})}{\tilde{g}_j(\tilde{x})},
\end{equation}
where $\tilde{f}_j,\tilde{g}_j\in{\cal O}_{K}[\tilde{x}]$ and 
$v\tilde{f}_j(\tilde{x})=0=v \tilde{g}_j (\tilde{x})$.
Consequently the first factor in the representation (\ref{represent}) is an 
element of ${\cal O}_{K}$ and we have shown that $z_1,\ldots,z_m\in 
{\cal O}_{K}[\tilde{x}]_q$ for the prime $q:={\cal O}_{K}[\tilde{x}]\cap
{\cal M}_{P}$.
\QED%
\begin{proposition}                         \label{isofstF}
Let $(F|K,P)$ be an immediate, valued function field of transcendence degree 
$1$ and assume that $F|K$ is separable and that $(K,P)$ is separably tame. 
Then $P$ is strongly smoothly ${\cal O}_{K}$-uniformizable.
\end{proposition}
\proof
By Theorem~\ref{stt3} there exists some $x\in F$ such that $(F,P)\subset 
(K(x)^h,P^h )$. Since $(K,P)$ is separably tame and hence separable-algebraically 
maximal, Lemma~\ref{close} shows that condition (\ref{trat}) holds in $(K(x)|K,P)$.
Therefore $P|_{K(x)}$ is strongly smoothly ${\cal O}_K$-uniformizable by 
Lemma~\ref{ist}. Lemma~\ref{etale+smoothly uniform}, (2) and (3) yield that 
$P$ is strongly smoothly ${\cal O}_{K(x)}$-uniformizable. The assertion now 
follows from Corollary \ref{trans}.
\QED
\subsection{Extensions within the completion}     
In this subsection the proof of the main result Theorem~\ref{MTdens} is 
provided. The subsequent two facts are main ingredients of that proof.
\begin{proposition}                               \label{comp1}
Let $(L|K,P)$ be a finitely generated, separable extension within the
completion of $(K,P)$. Then $P$ is strongly smoothly ${\cal O}_K$-uniformizable.
\end{proposition}
\proof
By assumption there exists a transcendence basis $T$ of $L|K$ such that 
$L|K(T)$ is separable-algebraic. By induction on the transcendence degree, 
using the Lemmata~\ref{close} and \ref{ist} and Corollary \ref{trans} we find 
that $P|_{K(T)}$ is strongly smoothly ${\cal O}_K$-uniformizable.

Since $L|K(T)$ is separable-algebraic and $L$ lies in the completion of $K$ 
which is also the completion of $K(T)$, $L$ must lie within the henselization 
of $K(T)$. Hence by Lemma~\ref{etale+smoothly uniform} $P$ is strongly smoothly
${\cal O}_{K(T)}$-uniformizable and thus again by Corollary \ref{trans}  
strongly smoothly ${\cal O}_K$-uniformizable.
\QED
\begin{lemma}			\label{imdlsep}
Every immediate extension of a defectless field is separable.
\end{lemma}
\proof
Let $(L|K,P)$ be an immediate extension and assume that $(K,P)$ is defectless. 
It suffices to show that every finite, purely inseparable extension $M$ of $K$ is 
linearly disjoint to $L$ over $K$. Let $e_M$ and $f_M$ be the ramification index 
and the residue degree of the unique extension of $P$ to $M$ and define similarly the 
ramification index $e_{L.M}$ and the residue degree $f_{L.M}$ in the extension $L.M|L$. 
The fundamental (in)equality then yields:
\[
[M:K]\geq [L.M:L]\geq e_{L.M}f_{L.M}\geq e_Mf_M =[M:K]
\]
and thus the assertion $[M:K]=[L.M:L]$.
\QED
\textbf{Proof of Theorem~\ref{MTdens}}: let $F_0$ be an intermediate field of 
$F|K$ such that $P|_{F_0}$ is an Abhyankar place and $(F,P)$ lies in the 
completion of $(F_0,P)$.

The valued field $(K,P)$ is defectless by assumption respectively because 
$P|_K=\mathrm{id}_K$. Hence the Generalized Stability Theorem [K7], Thm.~1 
yields that $(F_0,P)$ is defectless. The extension $F|F_0$ is immediate 
hence separable due to Lemma \ref{imdlsep}. Proposition \ref{comp1} thus 
yields that the place $P$ is strongly smoothly ${\cal O}_{F_0}$-uniformizable.

Let $Z\subset{\cal O}_{P}$ be a finite set. For every $z\in Z\cap{\cal M}_{P}$
we choose a representation $z=uz^\prime$ such that $u\in{\cal O}_{P}^\times$ 
and $z^\prime\in{\cal O}_{F_0}$ holds. Let $U\subset {\cal O}_{P}^\times$ and 
$Z^\prime\subset {\cal O}_{F_0}$ be the finite sets consisting of all of the 
elements $u$ and $z^\prime$ appearing in these representations. 
\smallskip

\noindent
\textbf{Case 1 of the theorem:} we apply Theorem~1.1 of [K--K] to obtain that 
$P|_{F_0}$ is strongly smoothly $K$-uniformizable. Corollary \ref{trans} then 
already yields that $P$ is strongly smoothly $K$-uniformizable. However to 
prove the existence of the morphism $f:X\rightarrow X_0$ we use case 1 of the 
Proposition \ref{stacking}: it yields the existence of a morphism 
$f: \Spec C\rightarrow\Spec B$ between affine $K$-models of $F$ and $F_0$ such 
that: 
\begin{itemize}
\item the $K$-algebra $B$ is smooth at $q_B:={\cal M}_{P}\cap B$ and 
$Z^\prime\subset B_{q_B}$,
\item $f$ is smooth at $q_C:={\cal M}_{P}\cap C$ and $U\subset C_{q_C}$.
\end{itemize}
Theorem~1.1 of [K--K] yields the existence of a regular system of parameters 
$(a_1,\ldots ,a_m)$ of $B_{q_B}$ such that every $z^\prime\in Z^\prime$ is a 
$B_{q_B}$-monomial in these parameters. Since the ring extension 
$C_{q_C}|B_{q_B}$ is flat, the elements $a_1,\ldots ,a_m$ remain a part of 
a regular system of parameters of $C_{q_C}$. Thus by construction every element 
$z=uz^\prime$ of $Z\subset C_{q_C}$ is a $C_{q_C}$-monomial in a regular system 
of parameters of $C_{q_C}$ as required. Using Theorem~1.1 of [K--K] a last time 
and the equation $vF_0=vF$ we get:
\[
\dim C_{q_C}\geq \dim B_{q_B}=\dim (vF_0\otimes\Q)=\dim (vF\otimes\Q).
\]
\noindent
\textbf{Case 2 of the theorem:} we apply Theorem~1.2 of [K--K] which yields 
that $P|_{F_0}$ is strongly $R$-uniformizable. Next we invoke case 2 of 
Proposition \ref{stacking} to obtain that $(P,Z)$ is $R$-uniformizable and 
the existence of a morphism $f: \Spec C\rightarrow\Spec B$ such that: 
\begin{itemize}
\item the $R$-algebra $B$ is regular at $q_B$ and $Z^\prime\subset B_{q_B}$,
\item the $B_{q_B}$-algebra $C_{q_B}$ is smooth and $U\subset C_{q_C}$.
\end{itemize}
The arguments used in case 1 to prove the assertions of the theorem carry 
over to case 2 just repacing Theorem~1.1 of [K--K] by Theorem~1.2. 
\QED
%
%
\section{Local uniformization by finite extension}           \label{sectlufe}
This section is devoted to the proofs of the main results Theorem 
\ref{main1}, \ref{main2} and \ref{main3}. In each of the three theorems local 
uniformization is achieved only after a finite extension of the function 
field in consideration. This finite extension can be choosen to be either 
Galois or an extension within a given separably tame extension of the 
function field. Although the proofs for the two cases are similar we present 
them separately to keep the exposition well-accessible. 
%
%
\subsection{The Proof of Theorem~\ref{main1}}           \label{sectluGal}
\textbf{Uniformization after a Galois extension}
\sn
We proceed by induction on the transcendence degree $n:=\trdeg E|K$ starting 
with the case $n=1$. Since by assumption $vE/vK$ is torsion and $EP|KP$ is 
algebraic Lemma~\ref{KsacP} implies that the extension $(E\sep|K\sep,{\cal P})$ 
and hence also its subextension $(E.K\sep|K\sep,{\cal P})$ are immediate. 
Since $(K\sep,{\cal P})$ is a separably tame field, we can apply 
Proposition~\ref{isofstF} to see that ${\cal P}|_{E.K\sep}$ is strongly 
smoothly ${\cal O}_{K\sep}$-uniformizable.

We express every $z\in Z$ in the form $z=uz^\prime$, $u\in 
{\cal O}_{E.K\sep}^{\times}$ and $z^\prime\in {\cal O}_{K\sep}$: let $U$
and $Z^\prime$ be the finite sets of elements $u$ and $z^\prime$ appearing
in these expressions. Moreover let $Z_g\subset{\cal O}_{P}$ be a 
finite set of generators of $E|K$. An application of Proposition 
\ref{base change} (3) yields the existence of a finite Galois extension 
${\cal K}|K$ with the following properties:
\begin{itemize}
\item $({\cal P}|_{\cal E},U\cup Z^\prime\cup Z_g)$ is smoothly 
${\cal O}_{\cal K}$-uniformizable, where ${\cal E}:=E.{\cal K}$,
\item ${\cal K}$ contains $K\sep\cap E$.
\end{itemize}
${\cal K}$ is algebraically closed in ${\cal E}$: $E|K$ is assumed to be 
separable, hence $K\sep\cap E$ is the algebraic closure of $K$ in $E$. We 
conclude that $E$ and $K\sep$ are linearly disjoint over $K\sep\cap E$, thus 
${\cal E}$ and $K\sep$ are linearly disjoint over ${\cal K}$, which yields 
the assertion.

Let $X$ be an ${\cal O}_{\cal K}$-model of ${\cal E}|{\cal K}$ that is  
smooth at the center $x$ of ${\cal P}$ and such that 
$U\cup Z^\prime\subset {\cal O}_{X,x}$ holds. Then $U\subset 
{\cal O}_{X,x}^\times$ and the factorizations $uz^\prime =z\in {\cal O}_{X,x}$ 
hold. Moreover $z^\prime\in {\cal O}_{X,x}\cap {\cal O}_{K\sep}=
{\cal O}_{\cal K}$, where the last equality holds because ${\cal K}$ is 
algebraically closed in ${\cal E}$.
 
Finally let $E_0|K$ be an arbitrary subextension of $E|K$ of transcendence 
degree $n-1=0$. Then $E_0|K$ is a finite separable extension, hence $E_0
\subseteq {\cal K}$ and the assertion is proved for $n=1$.
\smallskip

Let us now assume that $n>1$. We choose a subextension $E_0|K$ of $E|K$ of 
transcendence degree $n-1$ such that $E|E_0$ is separable. Such a subextension 
always exists: choose a separating transcendence basis $T$ of $E|K$ and a 
subset $T_0\subset T$ such that $\trdeg E|K(T_0)=1$. Set $E_0:=K(T_0)\subset 
E$, then $E|E_0$ is separable. 

Since $vE/vK$ is a torsion group and $EP|KP$ is algebraic, the same holds for 
$vE/vE_0$ and $EP|E_0P$. Hence by what we have already shown for the case $n=1$ 
and by the remarks on standard-smooth algebras following Theorem 
\ref{smoothness}, there exists a finite Galois extension ${\cal E}_0|E_0$ and 
an affine ${\cal O}_{{\cal E}_0}$-model $\Spec A$ of 
$E.{\cal E}_0 |{\cal E}_0$, $A\subset {\cal O}_{\cal P}$, with the following 
properties:
\begin{eqnarray}
\nonumber
\bullet & & A|{\cal O}_{{\cal E}_0}\mbox{ is standard-smooth,}\\
\label{factor1}
\bullet & & \forall z\in Z:\;\exists u\in A_{q_A}^{\times},
z^\prime\in{\cal O}_{{\cal E}_0}:\;\;z=uz^\prime,
\end{eqnarray}
where $q_A:=A\cap {\cal M}_{\cal P}$. Let $U\subset A_{q_A}^{\times}$ and 
$Z^\prime\subset {\cal O}_{{\cal E}_0}$ be the finite 
sets of elements $u$ and $z^\prime$ appearing in the factorizations 
(\ref{factor1}). 

Next we invoke Proposition \ref{descent 1} which yields a finitely generated 
${\cal O}_K$-algebra 
\begin{equation}
\label{universal descent}
S_0={\cal O}_K [x_1,\ldots ,x_r]\subseteq {\cal O}_{{\cal E}_0}
\end{equation} 
such that for every integrally closed domain 
$S^\prime\subseteq {\cal O}_{{\cal E}_0}$ with $S_0\subseteq S^\prime$ there 
exists a standard-smooth $S^\prime$-algebra $A^\prime\subseteq A_{q_A}$ with the 
property $U\subset (A^\prime)_{q^\prime}$, $q^\prime:=A^\prime\cap q_A$. 

Choose a valuation of $E\sep$ associated to ${\cal P}$; it is an extension of 
the valuation $v$ and we will denote it by $v$ too. As a direct consequence of 
the fact that $vE/vK$ is torsion and $EP|KP$ is algebraic we have that 
$v{\cal E}_0/vK$ is torsion and ${\cal E}_0{\cal P}|KP$ is algebraic. 

Let $E_1|K$ be a subextension of $E_0|K$ such that $E_0|E_1$ is separable and 
of transcendence degree $1$, then ${\cal E}_0|E_1$ is separable and of 
transcendence degree $1$ too. We apply the induction hypothesis to the valued 
function field $({\cal E}_0|K,{\cal P})$ and the subfunction field $E_1|K$: 
there exists a finite Galois extension ${\cal E}_1|E_1$ and a finite Galois extension 
${\cal K}|K$ within ${\cal E}_0.{\cal E}_1$ such that 
${\cal E}_0.{\cal E}_1|{\cal K}$ possesses an affine ${\cal O}_{\cal K}$-model 
$\Spec B$, $B\subset{\cal O}_{\cal P}$, with the following properties: 
\begin{eqnarray}
\label{smooth B}
\bullet & & B|{\cal O}_{\cal K}\mbox{ is smooth},\\
\label{universal generators}
\bullet & & \{x_1,\ldots ,x_r\}\subset B\mbox{ (see (\ref{universal descent}))},\\
\label{generators}
\bullet & & B\mbox{ contains a finite set of generators of }
{\cal E}_0|K,\\
\label{factor2}\bullet & & \forall z^\prime\in Z^\prime :\;\exists u^\prime\in 
B_{q_B}^\times , z^{\prime\prime}\in{\cal O}_{\cal K}:\;\; 
z^\prime =u^\prime z^{\prime\prime},
\end{eqnarray}
where $q_B:={\cal M}_{\cal P}\cap B$. Since $B_0:=B\cap{\cal E}_0$ is normal and 
contains $S_0$ (\ref{universal generators}) there exists a standard-smooth 
$B_0$-algebra $A_0\subseteq A_{q_A}$ with the property
\begin{equation}
\label{U}
U\subset (A_0)_{q_0}^\times,
\end{equation} 
where $q_0:=A_0\cap q_A$. The requirement (\ref{generators}) implies 
$\Frac B_0={\cal E}_0$ hence $\Frac A_0=\Frac A=E.{\cal E}_0$ by Proposition 
\ref{descent 1}.

We next consider the $B$-algebra $C:=A_0.B\subseteq {\cal O}_{\cal P}$: by 
Proposition \ref{ascent 1} it is standard-smooth, consequently $C$ is a smooth
${\cal O}_{\cal K}$-algebra due to (\ref{smooth B}). Furthermore we have $\Frac C=
E.{\cal E}_0.{\cal E}_1$ and since ${\cal E}_0$ and $E_0.{\cal E}_1$ are finite 
Galois extensions of $E_0\,$, so is ${\cal E}_0.{\cal E}_1\,$. 

Let $q_C:={\cal M}_{\cal P} \cap C$, the localization $C_{q_C}$ is a local 
extension of the ring $(A_0)_{q_0}$, hence $U\subset C_{q_C}^\times$ by (\ref{U}).
Similarly $C_{q_C}$ is a local extension of $B_{q_B}$ so that (\ref{factor2}) 
yields
\[
\forall z^\prime\in Z^\prime :\;\exists u^\prime\in 
C_{q_C}^\times , z^{\prime\prime}\in{\cal O}_{\cal K}:\;\; 
z^\prime =u^\prime z^{\prime\prime}.
\]
Combined with (\ref{factor1}) this shows that every $z\in Z$ can be factored 
in the form
\[
z=uu^\prime z^{\prime\prime}
\]
with $u,u^\prime\in C_{q_C}^\times$ and $z^{\prime\prime}\in{\cal O}_{\cal K}$.

Altogether we have shown that the ${\cal O}_{\cal K}$-model $\Spec C$ of 
$E.{\cal E}_0.{\cal E}_1 |{\cal K}$ fullfills the requirements stated in 
the assertion of Theorem \ref{main1}.
\QED
\bn
\textbf{Uniformization after an extension within a separably tame field}
\sn
The proof is similar to the one in the Galois case. We therefore put the focus 
on the differences between the two.
\smallskip

\noindent
We proceed by induction on the transcendence degree $n:=\trdeg E|K$ and start 
with the case $n=1$. Let $K'$ be the algebraic closure of $K$ within 
$E^{\rm st}$. Since by assumption $vE/vK$ is torsion and $EP|KP$ is 
algebraic Proposition~\ref{trac} implies that $(K',{\cal P})$ is a separably tame field 
and that the extension $(E^{\rm st}|K',{\cal P})$ and hence also its subextension 
$(E.K'|K',{\cal P})$ are immediate. Proposition~\ref{isofstF} now yields that 
${\cal P}|_{E.K'}$ is strongly smoothly ${\cal O}_{K'}$-uniformizable.

We define the finite sets $U\subset{\cal O}_{E.K'}^{\times}$, $Z^\prime\subset 
{\cal O}_{K'}$ and $Z_g\subset{\cal O}_{P}$ as in the proof of the Galois case. 
An application of Proposition \ref{base change} (2) yields the existence of a 
finite extension ${\cal K}|K$ within $K'$ such that:
\begin{itemize}
\item $({\cal P}|_{\cal E},U\cup Z^\prime\cup Z_g)$ is smoothly 
${\cal O}_{\cal K}$-uniformizable, where ${\cal E}:=E.{\cal K}$,
\item ${\cal K}$ is algebraically closed in ${\cal E}$.
\end{itemize}
The assertions of the theorem in the case $n=1$ now follow as in the Galois case.
\smallskip

Let us now assume that $n>1$ and that the assertion is true for transcendence 
degrees smaller than $n$. We take an arbitrary subextension $E_0|K$ of 
$E|K$ of transcendence degree $n-1$ such that $E|E_0$ is separable. 

As in the Galois case we deduce from what we have shown in the case $n=1$ 
the existence of a finite extension ${\cal E}_0|E_0$ within $E^{\rm st}$ such 
that $E.{\cal E}_0|{\cal E}_0$ possesses an affine 
${\cal O}_{{\cal E}_0}$-model $\Spec A$, $A\subset {\cal O}_{\cal P}$, with the 
properties
\begin{eqnarray}
\nonumber
\bullet & & A|{\cal O}_{{\cal E}_0}\mbox{ is standard-smooth,}\\
\label{factorst}
\bullet & & \forall z\in Z:\;\exists u\in A_{q_A}^{\times},
z^\prime\in{\cal O}_{{\cal E}_0}:\;\;z=uz^\prime,
\end{eqnarray}
where $q_A:=A\cap {\cal M}_{\cal P}$. Let $U$ and $Z^\prime$ be the finite 
sets of elements $u$ and $z^\prime$ appearing in the factorizations 
(\ref{factorst}).

Again choose a finitely generated ${\cal O}_K$-algebra  
\begin{equation}
\label{universal descent st}
S_0={\cal O}_K [x_1,\ldots ,x_r]\subseteq {\cal O}_{{\cal E}_0}
\end{equation} 
as described in Proposition \ref{descent 1}.

Let $E_1|K$ be a subsextension of $E_0|K$ such that $E_0|E_1$ is separable 
and of transcendence degree 1. Since the extension $E^{\rm st}|E$ is 
assumed to be separable-algebraic, the extension ${\cal E}_0|E_0$ is separable 
too. Hence ${\cal E}_0 |E_1$ is a separable extension of transcendence degree 1.

Let ${\cal E}_0^{\rm st}$ be the separable closure of ${\cal E}_0$ in 
$E^{\rm st}$. Since $E^{\rm st}{\cal P}|{\cal E}_0^{\rm st}{\cal P}$ is an 
algebraic extension, Lemma \ref{trac} shows that ${\cal E}_0^{\rm st}$ is a 
separably tame field.

We apply the induction hypothesis to the valued function field 
$({\cal E}_0|K,{\cal P})$, the subfunction field $E_1|K$ and the separably tame 
extension field ${\cal E}_0^{\rm st}\subseteq {\cal E}_0\sep$:
we obtain a finite extension ${\cal E}_1|E_1$ within ${\cal E}_0^{\rm st}$ and 
a finite extension ${\cal K}|K$ within ${\cal E}_0.{\cal E}_1\subseteq 
{\cal E}_0^{\rm st}$ such that ${\cal E}_0.{\cal E}_1|{\cal K}$ possesses an 
affine ${\cal O}_{\cal K}$-model $\Spec B$, $B\subset{\cal O}_{\cal P}$, with the 
following properties: 
\begin{eqnarray}
\label{smooth B st}
\bullet & & B|{\cal O}_{\cal K}\mbox{ is smooth},\\
\nonumber
\bullet & & \{x_1,\ldots ,x_r\}\subset B\mbox{ (see (\ref{universal descent st}))},\\
\nonumber
\bullet & & B\mbox{ contains a finite set of generators of }
{\cal E}_0|K,\\
\label{factor2 st}
\bullet & & \forall z^\prime\in Z^\prime :\;\exists u^\prime\in 
B_{q_B}^\times , z^{\prime\prime}\in{\cal O}_{\cal K}:\;\; 
z^\prime =u^\prime z^{\prime\prime},
\end{eqnarray}
where $q_B:={\cal M}_{\cal P}\cap B$. 

Next a standard-smooth $B$-algebra $C\subset {\cal O}_{\cal P}$ is constructed 
as in the Galois case. It satisfies $\Frac C=E.{\cal E}_0.{\cal E}_1$ and due 
to (\ref{smooth B st}) yields a smooth ${\cal O}_{\cal K}$-model $X$ of the 
function field $E.{\cal E}_0.{\cal E}_1 |{\cal K}$.

Using (\ref{factorst}) and (\ref{factor2 st}) as in the Galois case the 
required factorization of the elements $z\in Z$ in the local ring $C_{q_C}$, 
$q_C:={\cal M}_{\cal P}\cap C$, can be verified. 
\QED
%
%
\subsection{The proofs of the Theorems \ref{main2} and \ref{main3}}
We start by proving the existence of an intermediate field $F_0$ of $F|K$ 
as required in both Theorem \ref{main2} and Theorem \ref{main3}. We therefore do not 
assume that the place $P$ be trivial on $K$: 
by Proposition~\ref{tb} we can choose a separating transcendence basis of $F|K$, 
which contains elements $x_1,\ldots,x_{\rho},y_1,\ldots,y_{\tau}$ such that 
$\{v x_1+vK,\ldots,v x_{\rho}+vK\}$ is a maximal set of rationally independent 
elements in $vF/vK$, and $\{y_1P,\ldots, y_{\tau}P\}$ forms a
transcendence basis of $FP|KP$; let $F_0:= K(x_1,\ldots, x_{\rho},
y_1,\ldots,y_{\tau})$. The extension $F|F_0$ then is separable, $vF/vF_0$ 
is a torsion group and $FP|F_0 P$ is algebraic. Moreover $P|_{F_0}$ is an 
Abhyankar place by construction. 
\medskip

\noindent
\textbf{Proof of Theorem \ref{main2}}

\noindent
Let $F_0$ be an arbitrary intermediate field of $F|K$ with the 
properties required in Theorem \ref{main2}. Let ${\cal P}$ be an extension of $P$ 
to the separable closure $F\sep$ of $F$ and let $v$ be a valuation associated
to ${\cal P}$.

By Theorem~\ref{main1} there exists a finite extension ${\cal F}|F$ within 
$F\sep$ and a finite extension ${\cal F}_0|F_0$ within ${\cal F}$ such that 
the function field ${\cal F}|{\cal F}_0$ possesses an affine 
${\cal O}_{{\cal F}_0}$-model $\Spec A$ with the following properties:
\begin{eqnarray}
\nonumber
\bullet & & A\mbox{ is smooth at }q_A:=A\cap {\cal M}_{\cal P},\\
\label{factor up}
\bullet & & \forall z\in Z:\;\exists u\in A_{q_A}^\times ,
z^{\prime}\in {\cal O}_{{\cal F}_0}:\;\; z=uz^\prime .
\end{eqnarray}
Let $U\subset A_{q_A}^\times$ and $Z^\prime\subset {\cal O}_{{\cal F}_0}$ 
be the finite sets of elements appearing in these factorizations.

Next we choose a finitely generated $K$-algebra
\begin{equation}
\label{K-algebra}
S_0=K[x_1,\ldots ,x_r]\subseteq {\cal O}_{{\cal F}_0}
\end{equation}
according to Proposition \ref{descent 1} applied to the 
${\cal O}_{{\cal F}_0}$-algebra $A$ and the finite set $U$.

By the choice of $F_0$ the place ${\cal P}|_{{\cal F}_0}$ 
is an Abhyankar place of ${\cal F}_0|K$. In particular the extension 
${\cal F}_0{\cal P}|K$ is finitely generated. Thus there exists a 
finite purely inseparable extension ${\cal K}|K$ such that  
the extension ${\cal F}_0{\cal P}.{\cal K}|{\cal K}$ is separable. 
Consequently, replacing $K$ by ${\cal K}$, $F_0$ by $F_0.{\cal K}$ and 
$F$ by $F.{\cal K}$ we can assume that already the residue field 
extension ${\cal F}_0{\cal P}|K$ is separable.

We can thus apply 
Theorem~1.1 of [K--K] to the valued function field $({\cal F}_0|K,{\cal P})$: 
there exists an affine, smooth $K$-model $X_0=\Spec B$ of ${\cal F}_0|K$, 
$B\subset {\cal O}_{{\cal F}_0}$, and a regular parameter system 
$(a_1,\ldots ,a_d)$ of $B_{q_B}$, $q_B:={\cal M}_{{\cal F}_0}\cap B$, with 
the properties:
\begin{eqnarray}
\label{constants}
\bullet & & Z^\prime\cup\{x_1,\ldots ,x_r\}\subset B\;
\mbox{ (see (\ref{K-algebra}))},\\
\label{monomials}
\bullet & & \mbox{every $z^\prime\in Z^\prime$ is a 
$B_{q_B}$-monomial in $\{a_1,\ldots ,a_d\}$},\\
\label{dimensions}
\bullet & & \dim B_{q_B}=\dim (v{\cal F}_0 \otimes\Q ). 
\end{eqnarray}
Property (\ref{constants}) implies $S_0\subseteq B$, therefore due to 
Proposition \ref{descent 1} and since $B$ is normal, there exists a finitely 
presented $B$-algebra $C\subseteq A_{q_A}$ with the properties:
\begin{eqnarray}
\nonumber
\bullet & & C\mbox{ is smooth at }q_C:=C\cap {\cal M}_{\cal P} 
\mbox{ and }\Frac C={\cal F},\\ 
\label{U in C}
\bullet & & U\subset C_{q_C}.
\end{eqnarray}
Note that since $A_{q_A}|C_{q_C}$ is a local extension (\ref{U in C}) and 
the definition of $U$ yield $U\subset C_{q_C}^\times$. Consequently 
(\ref{factor up}) and (\ref{monomials}) show that every $z\in Z$ is a 
$C_{q_C}$-monomial in $\{a_1,\ldots ,a_d\}$. Since the extension $C|B$ 
is smooth the local extension $C_{q_C}|B_{q_B}$ is flat, hence 
$\{a_1,\ldots ,a_d\}$ is part of a regular parameter system of $C_{q_C}$, 
[M], Thm.~23.7. Moreover using (\ref{dimensions}) we get
\[
\dim C_{q_C}\geq \dim B_{q_B}=\dim (v {\cal F}_0\otimes\Q )=
\dim (v F\otimes\Q ),
\]
where the last equality holds because $v{\cal F}_0/vF_0$ and $vF/vF_0$ are 
torsion groups.

Altogether we have shown that the smooth $K$-morphism $f:\Spec C=:X\rightarrow 
X_0$ induced by the extension $C|B$ fullfills the requirements stated in the 
assertions of Theorem \ref{main2}.
\QED
\bn
\textbf{Proof of Theorem~\ref{main3}}
\sn
The proof is very similar to that of Theorem~\ref{main2}. Using the same notation 
we therefore only point out the differences between the two.
\begin{itemize}
\item In (\ref{K-algebra}) $S_0$ is choosen to be a finitely generated $R$-algebra 
$R[x_1,\dots ,x_r]\subseteq {\cal O}_{{\cal F}_0}$.
\item Theorem 1.2 of [K--K] is used to obtain an affine $R$-model $X_0=\Spec B$ 
of the function field ${\cal F}_0|K$ that is regular at the center $q_B$ of 
${\cal P}$ on $X_0$. The requirements for an application of this theorem are 
satisfied by assumption except for the separability of ${\cal F}_0{\cal P}|KP$ 
that follows from the assumed perfectness of $KP$. 
\item The dimension formula (\ref{dimensions}) has to be replaced with
\[ 
\dim B_{q_B}=\left\{
\begin{array}{rl}
\dim (v{\cal F}_0/vK\otimes\Q)+1 & \mbox{ if $\dim R=1$ or $\trdeg (KP|R/M)>0$}\\
\dim (v{\cal F}_0/vK\otimes\Q)+2 & \mbox{ in the remaining cases}
\end{array}\right.
\]
\item Note that $B$ is not necessarily normal so that in the construction of the 
$B$-algebra $C$ we have to add an intermediate step: using Proposition 
\ref{descent 1} and (\ref{constants}) we obtain a finitely presented 
$B_{q_B}$-algebra $C^\prime\subseteq A_{q_A}$ that is smooth at $q_{C^\prime}=
{\cal M}_P\cap C^\prime$ and satisfies $U\subset C^\prime_{q_{C^\prime}}$. For 
$C^\prime =B_{q_B}[x_1,\ldots ,x_r]$ we then set $C :=B[x_1,\ldots ,x_r]$ and 
get $C_{q_B}=C^\prime$ and $C_{q_C}=C^\prime_{q_{C^\prime}}$. The smoothness 
of $C^\prime |B_{q_B}$ at $q_{C^\prime}$ and the regularity of $B_{q_B}$ imply 
the regularity of $C_{q_C}$
\end{itemize}
Altogether it is shown that the $R$-morphism $f:\Spec C=:X\rightarrow 
X_0:=\Spec B$ induced by the extension $C|B$ fullfills the requirements stated 
in the assertions of Theorem \ref{main3}.
\QED
\medskip

\noindent
{\bf References}
\newenvironment{reference}%
{\begin{list}{}{\setlength{\labelwidth}{5em}\setlength{\labelsep}{0em}%
\setlength{\leftmargin}{5em}\setlength{\itemsep}{-1pt}%
\setlength{\baselineskip}{3pt}}}%
{\end{list}}
\newcommand{\lit}[1]{\item[{#1}\hfill]}
\begin{reference}
\lit{[A1]} {Abhyankar, S.$\,$: {\it Local uniformization on algebraic
surfaces over ground fields of characteristic $p\neq 0$}, Ann.\ Math.
\textbf{63} (1956), 491--526}
\lit{[A2]} {Abhyankar, S.$\,$: {\it Resolution of singularities of
arithmetical surfaces}, Arithmetical Algebraic Geometry, Harper and Row,
New York (1965)}
\lit{[B]} {Bourbaki, N.$\,$: {\it Commutative algebra}, Paris (1972)}
\lit{[dJ]} {de Jong, A.~J.$\,$: {\it Smoothness, semi-stability and
alterations}, Publ. Math. IHES {\bf 83} (1996), 51--93}
\lit{[EGA IV]} {Grothendieck, A.: {\it \'Etude locale de sch\'emas 
et des morphismes de sch\'emas}, Publ. Math. IHES {\bf 20} (1964), 
{\bf 24} (1965), {\bf 28} (1966), {\bf 32} (1967)}
\lit{[EL]} {Elliott, G.~A.$\,$: {\it On totally ordered groups, and
$K_0$}, in: Ring Theory Waterloo 1978, eds.\ D.~Handelman and
J.~Lawrence, Lecture Notes Math.\ {\bf 734}, 1--49}
\lit{[EN]} {Endler, O.$\,$: {\it Valuation theory}, Berlin (1972)}
\lit{[J--R]} {Jarden, M.\ -- Roquette, P.$\,$: {\it The Nullstellensatz
over $\wp$--adically closed fields}, J.\ Math.\ Soc.\ Japan {\bf 32}
(1980), 425--460}
\lit{[KA]} {Kaplansky, I.$\,$: {\it Maximal fields with valuations I},
Duke Math.\ J.\ {\bf 9} (1942), 303--321}
\lit{[K--K]} {Knaf, H.\ -- Kuhlmann, F.--V.$\,$: {\it Abhyankar places
admit local uniformization in any characteristic}, to appear}
\lit{[K1]} {Kuhlmann, F.--V.$\,$: {\it On local uniformization in
arbitrary characteristic}, The Fields Institute Preprint Series, Toronto
(1997)}
\lit{[K2]} {Kuhlmann, F.--V.$\,$: {\it Valuation theoretic and model
theoretic aspects of local uniformization}, in:
Resolution of Singularities - A Research Textbook in Tribute to Oscar
Zariski. Herwig Hauser, Joseph Lipman, Frans Oort, Adolfo Quiros
(eds.), Progress in Mathematics Vol.\ {\bf 181}, Birkh\"auser
Verlag Basel (2000), 381--456}
\lit{[K3]} {Kuhlmann, F.-V.: {\it Elementary properties of power series
fields over finite fields}, J.\ Symb.\ Logic {\bf 66} (2001), 771--791}
\lit{[K4]} {Kuhlmann, F.-V.: {\it Value groups, residue fields and bad
places of rational function fields}, to appear in: Trans.\ Amer.\ Math.\
Soc.}
\lit{[K5]} {Kuhlmann, F.--V.$\,$: {\it Places of algebraic function
fields in arbitrary characteristic}, Adv. Math. {\bf 188} (2004), 399--424}
\lit{[K6]} {Kuhlmann, F.--V.$\,$: {\it A classification of Artin
Schreier defect extensions and a characterization of defectless fields},
in preparation}
\lit{[K7]} {Kuhlmann, F.--V.$\,$: {\it Elimination of Ramification in
Valued Function Fields I: The Generalized Stability Theorem}, in
preparation}
\lit{[K8]} {Kuhlmann, F.--V.$\,$: {\it Elimination of Ramification in
Valued Function Fields II: Henselian Rationality}, in preparation}
\lit{[K--P--R]} {Kuhlmann, F.-V.\ -- Pank, M.\ -- Roquette, P.$\,$:
{\it Immediate and purely wild extensions of valued fields},
manuscripta math.\ {\bf 55} (1986), 39--67}
\lit{[Kun]} {Kunz, E.: {\it K\"ahler Differentials}, Braunschweig (1986)}
\lit{[L]} {Lang, S.$\,$: {\it Algebra}, New York (1965)}
\lit{[M]} {Matsumura, H.$\,$: {\it Commutative Ring Theory}, Cambridge (1986)}
\lit{[Ray]} {Raynaud, M.$\,$: {\it Anneaux Locaux Hens\'eliens}, 
Berlin--Heidelberg--New York (1970)}
\lit{[R]} {Ribenboim, P.$\,$: {\it Th\'eorie des valuations}, Les
Presses de l'Uni\-versit\'e de Mont\-r\'eal (1964)}
\lit{[S]} {Spivakovsky, M.$\,$: {\it Resolution of singularities I:
local uniformization}, manu\-script, Toronto (1996)}
\lit{[W]} {Warner, S.$\,$: {\it Topological fields}, Mathematics
studies {\bf 157}, North Holland, Amsterdam (1989)}
\lit{[Z1]} {Zariski, O.$\,$: {\it Local uniformization on
algebraic varieties}, Ann.\ Math.\ {\bf 41} (1940), 852--896}
\lit{[Z2]} {Zariski, O.$\,$: {\it The reduction of singularities of an
algebraic surface}, Ann.\ Math.\ {\bf 40} (1939), 639--689}
\lit{[Z3]} {Zariski, O.$\,$: {\it A simplified proof for resolution of
singularities of an algebraic surface}, Ann.\ Math.\ {\bf 43} (1942),
583--593}
\lit{[Z--S]} {Zariski, O.\ -- Samuel, P.$\,$: {\it Commutative
Algebra}, Vol.\ II, New York--Heidel\-berg--Berlin (1960)}
\end{reference}
\adresse
\end{document}